\documentclass[journal]{IEEEtran}

\usepackage{cite}
\usepackage{amsmath,amssymb,amsfonts,bm}
\usepackage{bbm}
\usepackage{graphicx}
\usepackage{subfigure}
\usepackage{textcomp}
\usepackage{algorithm}
\usepackage{algpseudocode}

\makeatletter
\renewcommand*\env@matrix[1][*\c@MaxMatrixCols c]{%
	\hskip -\arraycolsep
	\let\@ifnextchar\new@ifnextchar
	\array{#1}}
\makeatother
\usepackage{hyperref}
\usepackage{enumitem}
\newcommand{\tabincell}[2]{\begin{tabular}{@{}#1@{}}#2\end{tabular}}
\def \qed {\hfill \vrule height6pt width 6pt depth 0pt}
\usepackage{mathrsfs}
\usepackage{xcolor}
\allowdisplaybreaks[4]
\newtheorem{definition}{Definition}
\newtheorem{theorem}{Theorem}
\newtheorem{lemma}{Lemma}
\newtheorem{assum}{Assumption}

\newtheorem{remark}{Remark}

\newtheorem{problem}{Problem}
\newenvironment{proof}{\begin{IEEEproof}}{\end{IEEEproof}}

\begin{document}
	\title{Policy Gradient Adaptive Control for the LQR: Indirect and Direct Approaches 
		\thanks{This work
			was supported by ETH Zurich and the SNF through the NCCR Automation.}
		\thanks{F. Zhao and F. D\"{o}rfler are with the Department of Information Technology and Electrical Engineering, ETH Z\"{u}rich, 8092 Z\"{u}rich, Switzerland. (e-mail: zhaofe@control.ee.ethz.ch; dorfler@control.ee.ethz.ch)}\thanks{A. Chiuso is with the Department of Information Engineering, University of Padova, Via Gradenigo 6/b, 35131 Padova, Italy. (e-mail: alessandro.chiuso@unipd.it)}
		\author{Feiran Zhao, Alessandro Chiuso, Florian D\"{o}rfler}}
	\maketitle
	
	\begin{abstract}
	Motivated by recent advances of reinforcement learning and direct data-driven control, we propose policy gradient adaptive control (PGAC) for the linear quadratic regulator (LQR), which uses online closed-loop data to improve the control policy while maintaining stability. Our method adaptively updates the policy in feedback by descending the gradient  of the LQR cost and is categorized as indirect, when gradients are computed via an estimated model, versus direct, when gradients are derived from data using sample covariance parameterization. Beyond the vanilla gradient, we also showcase the merits of the natural gradient and Gauss-Newton methods for the policy update. Notably, natural gradient descent bridges the indirect and direct PGAC, and the Gauss-Newton method of the indirect PGAC leads to an adaptive version of the celebrated Hewer's algorithm. To account for the uncertainty from noise, we propose a regularization method for both indirect and direct PGAC. For all the considered PGAC approaches, we show  closed-loop stability and convergence of the policy to the optimal LQR gain. Simulations validate our theoretical findings and demonstrate the robustness and computational efficiency of PGAC.
	\end{abstract}
	
	
	\section{Introduction}
	The history of adaptive control is almost as long as the entire control field \cite{annaswamy2021historical}, and it is currently revived in the context of reinforcement learning (RL)~\cite{annaswamy2023adaptive}. A fundamental principle of an adaptive control system is its ability to monitor its own performance and adjust its parameters in the direction of better performance~\cite{drenick1957adaptive}. Adaptive control for linear time-invariant systems with unknown parameters is widely studied, and the manifold approaches can be divided with three orthogonal classifications. A commonly adopted one is indirect (when a dynamical model is identified followed by model-based control), versus direct (when bypassing identification)~\cite{ioannou1996robust}. Based on control objectives, approaches are categorized as seeking stability (i.e., convergence of signals) or optimality (i.e., convergence of a performance index). Another perspective considers the policy update rule: one-shot-based methods solve an online optimization problem to obtain the policy, whereas gradient-based methods update the policy iteratively using online gradient information.

	
	Classical adaptive control methods mainly seek robust stability, including indirect and direct model-reference adaptive control (MRAC), self-tuning regulators, and model-free adaptive control~\cite{ioannou1996robust,8621060}. In particular, MRAC designs the control input based on Lyapunov theory to guarantee stability. In comparison, other approaches focus on optimality, aligning more closely with the Zames's definition of adaptive control \cite{zames1998adaptive}, i.e., improve over the best performance with \textit{prior information}. A widely adopted  optimality criterion is the infinite-horizon linear quadratic regulator (LQR) cost, which quantifies the $\mathcal{H}_2$ norm of the closed-loop system~\cite{lewis2009reinforcement, bian2016value,mania_certainty_2019,lu2023almost,simchowitz2020naive,cohen2019learning,chekan2024fully,karafyllis2019adaptive}. A representative instance is adaptive dynamical programming (ADP), a classical approach of reinforcement learning (RL)~\cite{lewis2009reinforcement, bian2016value}. While it shows convergence to the optimal LQR gain, the stability of the closed-loop system remains unclear. An exception is the identification-based policy iteration method~\cite{song2025robustness}, which shows stability under sufficiently small identification error.  Recently, there have been adaptive control methods showing both stability and convergence of the policy, based on the indirect certainty-equivalence LQR~\cite{mania_certainty_2019,lu2023almost,simchowitz2020naive,cohen2019learning,chekan2024fully,karafyllis2019adaptive}. 
	These optimality-seeking methods commonly require persistently exciting (PE) inputs, without which the bursting phenomenon may happen due to insufficient excitation~\cite{anderson2005failures}.

	
	A representative example of one-shot-based methods is indirect adaptive LQR, where the state-feedback policy is determined as the Riccati or semi-definite programming (SDP) solution of a certainty-equivalence LQR problem~\cite{mania_certainty_2019, karafyllis2019adaptive, lu2023almost,simchowitz2020naive,cohen2019learning,chekan2024fully}. 	
	However, one-shot-based methods are sensitive to uncertainty, as the policy may significantly vary with large noise. Gradient-based methods are widely adopted in classical adaptive control, which uses different objectives for online gradient computation  \cite{annaswamy2021historical}. For example, the well-known MIT-rule optimizes the squared one-step prediction error, and MRAC includes an additional objective that reflects dynamics of the system. Compared with one-shot-based methods, gradient-based methods are computationally more efficient, as only a single gradient must be computed instead of a complete optimization problem. Moreover, they are robust to noise due to smooth and incremental policy updates. However, the existing gradient-based approaches mainly focus on stability, and it remains unclear if the LQR objective can be adopted to ensure optimality in adaptive control~\cite{annaswamy2023adaptive}. A major challenge is the non-convexity of the LQR cost in the state-feedback gain,  hindering proving convergence of gradient-based methods.
	
	
	Theoretical studies of the policy gradient method, a foundational cornerstone of modern RL~\cite{bin2022towards}, have shed a light on the non-convex optimization landscape of the LQR~\cite{fazel2018global,mohammadi2022convergence,zhao23global,bin2022towards}. Specifically, the LQR cost is known to be \textit{gradient dominated} in the state-feedback gain, leading to global linear convergence for policy gradient methods~\cite{fazel2018global}. However, computing the policy gradient requires an exact dynamical model. While zeroth-order methods can be used for gradient estimates in the absence of the model, they are intrinsically unsuitable for adaptive control, as the gradient can be estimated only \textit{after} observing multiple long trajectories~\cite{mohammadi2022convergence,zhao2024data}. Even with an exact policy gradient, the \textit{sequential stability} of the closed-loop system under \textit{switching policies} remains unclear. In fact, proving the stability of control systems based on RL (e.g., ADP and policy gradient methods)
	is a difficult task \cite{annaswamy2023adaptive}, and the intersection of RL and adaptive control continues to be a compelling and largely unexplored direction~\cite{annaswamy2023integration}.

	\begin{figure}[t]
		\centerline{\includegraphics[width=70mm]{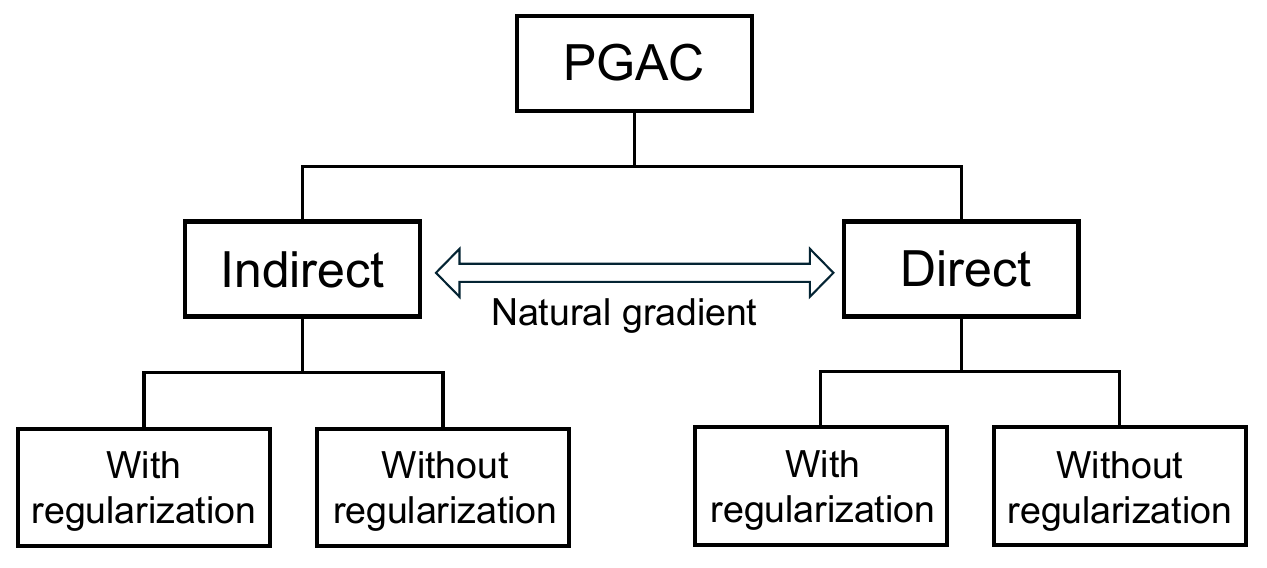}}
		\caption{Classes of policy gradient adaptive control.}
		\label{pic:overview}
	\end{figure}
	
	Recently, there has been a growing interest in direct data-driven control inspired by subspace methods and behavioral system theory \cite{markovsky2021behavioral, coulson2019data}. With a data-based system parameterization, direct methods obtain the LQR gain directly from a single batch of PE data~\cite{de2019formulas}. Regularization methods are introduced to promote certainty-equivalence or robustness~\cite{dorfler2021certainty,dorfler22on} for direct LQR design or to harness the uncertainty for the separation principle of predictive control \cite{chiuso2023harnessing}.  Despite these advances, the adaptation of direct methods is acknowledged as a fundamental open problem in the comprehensive surveys~\cite{markovsky2021behavioral,annaswamy2023adaptive}. By proposing a sample covariance parameterization for the LQR and developing a \textit{projected gradient dominance} property, our previous works~\cite{zhao2023data,zhao2024data} attack the direct adaptive control problem with a gradient-based method called data-enabled policy optimization (DeePO), which has seen successful applications in power converter systems~\cite{zhao2024direct}, aerospace control~\cite{10981973}, and autonomous bicycle control~\cite{persson2025adaptive}. However, the closed-loop stability of DeePO remains unclear.
	
	In this article, we propose policy gradient adaptive control (PGAC) for the LQR, which is a unified gradient-based framework touching upon both indirect and direct approaches and seeking  both stability and optimality. Starting from a batch of offline data and an initial stabilizing gain, it alternates between control, where the input signifies a state-feedback term plus probing noise ensuring persistency of excitation, and gradient-based policy update, where the policy gradient of the LQR cost is approximately computed with online closed-loop data. Our prior DeePO method~\cite{zhao2024data} is a special instance: namely,  direct PGAC, using the vanilla gradient of the certainty-equivalent covariance-parameterized LQR for policy update. For indirect PGAC, we use the certainty-equivalence LQR cost with ordinary least-squares identification to compute the policy gradient. Beyond the vanilla policy gradient, we develop natural gradient and Gauss-Newton methods for PGAC. Notably, natural gradient descent bridges   indirect and direct PGAC, and the Gauss-Newton method yields an adaptive variant of the celebrated Hewer's algorithm \cite{hewer1971iterative},  an iterative method to solve the Riccati equation and a cornerstone of ADP. Our adaptive Hewer's algorithm coincides with the online identification-based policy iteration algorithm~\cite{song2025robustness}. Since the certainty-equivalence LQR formulation disregards uncertainty in the estimated model and may lead to an unstable closed-loop system \cite{chiuso2023harnessing, zhao2025regularization}, we propose a variance-based regularization method for both indirect and direct PGAC to account for the uncertainty. The classification of our certified PGAC methods is shown in Fig. \ref{pic:overview}.
	
	Motivated by Zames’s perspective that adaptation and learning involve the acquisition of \textit{information} about the plant~\cite{zames1998adaptive}, we adopt the \textit{signal-to-noise ratio} (SNR) of the collected data as a quantitative information metric, reflecting the uncertainty in the identified dynamics. This notion satisfies \textit{Zames’s first monotonicity principle}~\cite{zames1998adaptive}, as the SNR typically increases monotonically over time. For all the aforementioned PGAC approaches, we show that the optimality gap of the LQR gain is upper bounded by two terms signifying an exponential decrease in the initial optimality gap plus a bias scaling linearly with inverse of the SNR. This aligns with \textit{Zames’s second monotonicity principle} of adaptive control \cite{zames1998adaptive}, i.e., the performance improves with increasing information. Our results capture the convergence rate of the policy as a function of the SNR, improving upon the rate established in DeePO \cite{zhao2024data}. For example, the optimality gap decreases with time $t$ as $\mathcal{O}(1/\sqrt{t})$ for $\text{SNR}_t \sim \mathcal{O}(\sqrt{t})$, which corresponds to the case of Gaussian noise and a constant excitation level. 
	
	
	Apart from non-asymptotic policy convergence certificates, we provide closed-loop stability guarantees of PGAC by leveraging a favorable feature of gradient-based methods, i.e., the update rate can be easily controlled via the stepsize. Specifically, we show that under a sufficiently small stepsize, the closed-loop system is sequentially stable~\cite{cohen2019learning} under switching polices, and the state is upper-bounded by two terms, i.e., an exponential decrease in the initial state plus an upper bound of probing and process noise. While the one-shot-based adaptive control methods also ensure stability, they require a dwell time to mitigate the burn-in effect of switching policies and a safety mechanism in case of divergence~\cite{ chekan2024fully, cohen2019learning, lu2023almost}. We believe that our analysis techniques provide a viable solution to proving closed-loop stability in RL algorithms \cite{annaswamy2023adaptive}.
	
	We perform simulations on a benchmark problem in \cite{dean2020sample} to validate our theoretical findings and demonstrate favorable computational efficiency and robustness of PGAC. 
	 
	The remainder of this article is organized as follows. Section \ref{sec:prob} recapitulates data-driven formulations of the LQR. Section \ref{sec:PGAC} formulates the adaptive control problem and proposes PGAC. Section \ref{sec:grad} proposes two variants of gradient descent for PGAC. Section \ref{sec:reg} develops a regularization method for PGAC. Section \ref{sec:simu} uses simulations to validate our theoretical results. Concluding remarks are made in Section \ref{sec:conclu}. All proofs are provided in the Appendix.




	\textbf{Notation.} We let $I_n$ denote the $n$-by-$n$ identity matrix. We let $\underline{\sigma}(\cdot)$ and $\sigma_n(\cdot)$  denote the minimal and $n$-th large singular value of a matrix, respectively. We let $\|\cdot\|$  denote the $2$-norm of a vector or matrix, and $\|\cdot\|_F$ the Frobenius norm. We let $A^\dagger:=A^{\top}(AA^{\top})^{-1}$ denote the right inverse of a full row rank matrix $A\in \mathbb{R}^{n\times m}$. We let $\mathcal{N}(A)$  denote the nullspace of $A$ and $\Pi_A := I_{m}-A^{\dagger}A$ the projection operator onto $\mathcal{N}(A)$.
		
	\section{Data-driven formulations of the LQR}\label{sec:prob}
This section recapitulates indirect certainty-equivalence LQR with least-squares identification \cite{dorfler22on} and direct certainty-equivalence LQR   with   covariance parameterization \cite{zhao2024data}.
	\subsection{Model-based LQR}
 
	Consider a discrete-time linear system
	\begin{equation}\label{equ:sys}
	\left\{\begin{aligned}
	x_{t+1} & =A x_t+B u_t+w_t \\
	z_t & =\begin{bmatrix}
	Q^{1 / 2} & 0 \\
	0 & R^{1 / 2}
	\end{bmatrix}
	\begin{bmatrix}
	x_t \\
	u_t
	\end{bmatrix}
	\end{aligned}\right.,
	\end{equation}
	where $t\in \mathbb{N}$, $x_t\in\mathbb{R}^{n}$ is the state, $u_t\in\mathbb{R}^{m}$ is the control input, $w_t \in \mathbb{R}^n$ is the noise, and $z_t$ is the performance signal of interest. We assume that $(A,B)$ are controllable and the weighting matrices $(Q, R)$ are positive definite.  
	
	The LQR problem is phrased as finding a state-feedback gain $K\in \mathbb{R}^{m\times n}$ that minimizes the $\mathcal{H}_2$-norm of the transfer function $\mathscr{T}(K):w \rightarrow z$ of the closed-loop system
	\begin{equation}\label{equ:closedsys}
	\begin{bmatrix}
	x_{t+1} \\
	z_t
	\end{bmatrix}=\begin{bmatrix}[c|c]
	A+BK & I_n \\
	\hline \begin{bmatrix}
	Q^{1 / 2} \\
	R^{1 / 2} K
	\end{bmatrix} & 0
	\end{bmatrix}\begin{bmatrix}
	x_t \\
	w_t
	\end{bmatrix}.
	\end{equation}
	When $A+BK$ is stable, it holds that \cite{anderson2007optimal}
	\begin{equation}\label{equ:transfer}
	\|\mathscr{T}(K)\|_2^2  = \text{Tr}((Q+K^{\top}RK)\Sigma)=:C(K),
	\end{equation}
	where $\Sigma$ is the closed-loop state covariance matrix obtained as the positive definite solution to the Lyapunov equation
	\begin{equation}\label{equ:Sigma}
	\Sigma = I_n + \left(A+BK\right)\Sigma \left(A+BK\right)^{\top}.
	\end{equation}
	We refer to $C(K)$ as the LQR cost and to (\ref{equ:transfer})-(\ref{equ:Sigma}) as a \textit{policy parameterization} of the LQR.

	For known $(A,B)$, there are alternative formulations to find the optimal LQR gain $K^*:=\arg\min_{K}C(K)$ of (\ref{equ:transfer})-(\ref{equ:Sigma}), e.g., via the celebrated Riccati equation \cite{anderson2007optimal}. We recall a well-known iterative algorithm to solve the Riccati equation, called \textit{Hewer's algorithm}~\cite{hewer1971iterative} (also known as policy iteration in reinforcement learning (RL) \cite{lewis2009reinforcement}), which starts from an initial stabilizing policy $K_0$ and alternates between
	\begin{subequations}\label{equ:hewer}
		\begin{align}
			&P_{i+1} = Q + K_i^{\top}RK_i +(A+ BK_i)^{\top}P_{i+1}(A + BK_i), \label{equ:hewer1}\\
			&K_{i+1} = (R + B^{\top}P_{i+1}B)^{-1}B^{\top}P_{i+1}A, \label{equ:hewer2}
		\end{align} 
	\end{subequations}
	where \eqref{equ:hewer1} is the policy evaluation, and \eqref{equ:hewer2} is the policy improvement. The Hewer's algorithm has asymptotic convergence guarantees to the solution of the Riccati equation and the optimal policy $K^*$, which are fixed points of \eqref{equ:hewer}~\cite{hewer1971iterative}.
	
	For unknown $(A,B)$, there is a plethora of data-driven methods to find $K^*$, some of which we recapitulate below.
	
	
	%

	\subsection{Indirect certainty-equivalence LQR with ordinary least-square identification}
	The conventional approach to data-driven LQR design follows the certainty-equivalence principle: it first identifies nominal system matrices $(A,B)$ from data, and then solves the LQR problem regarding the identified model as the ground-truth. The identification step is based on subspace relations among the state-space data. Consider a $t$-long time series of states, inputs, unknown noises, and successor states
	\begin{equation}\label{equ:data}
	\begin{aligned}
	X_{0} &:= \begin{bmatrix}
	x_0& x_1& \dots& x_{t-1}
	\end{bmatrix}\in \mathbb{R}^{n\times t},\\
	U_{0} &:= \begin{bmatrix}
	u_0& u_1& \dots& u_{t-1}
	\end{bmatrix}\in \mathbb{R}^{m\times t}, \\
	W_{0} &:= \begin{bmatrix}
	w_0& w_1& \dots& w_{t-1}
	\end{bmatrix}\in \mathbb{R}^{n\times t}, \\
	X_{1} &:= \begin{bmatrix}
	x_1& x_2& \dots& x_t
	\end{bmatrix}\in \mathbb{R}^{n\times t},
	\end{aligned}
	\end{equation}
	which satisfy the system dynamics
	\begin{equation}\label{equ:dynamics}
	X_1 = AX_0+ BU_0 + W_0.
	\end{equation}
	
	We assume that the data is {\em persistently exciting (PE)} \cite{willems2005note}, i.e., the block matrix of input and state data 
	$$
	D_0 := \begin{bmatrix}
	U_0 \\
	X_0
	\end{bmatrix}\in \mathbb{R}^{(m+n)\times t}
	$$
	has full row rank
	\begin{equation}\label{equ:rank}
	\text{rank}(D_0) = m+n,
	\end{equation}
	which is necessary for data-driven LQR design~\cite{van2023informativity}.
	
	Based on the subspace relations \eqref{equ:dynamics} and the rank condition (\ref{equ:rank}), an estimated model $(\hat{A},\hat{B})$  can be obtained as the unique solution to the ordinary least-squares problem
	\begin{equation}\label{equ:sysid}
	[\hat{B}, \hat{A}] = \underset{B, A}{\arg \min }\left\|X_1-[B,A] D_0\right\|_F = X_1D_0^{\dagger}.
	\end{equation}
	Following the certainty-equivalence principle~\cite{dorfler2021certainty}, the system $(A,B)$ is replaced with its estimate $(\hat{A},\hat{B})$ in (\ref{equ:transfer})-(\ref{equ:Sigma}), and the LQR problem can be reformulated as 
	\begin{equation}\label{prob:indirect} 
	\mathop{\text{minimize}}\limits_{K} ~ \text{Tr}\left((Q+K^{\top}RK)\Sigma\right), 
	\end{equation}
	where $\Sigma$ is the positive definite solution to
	\begin{equation}\label{equ:indirect_sigma}
		\Sigma = I_n + \left(\hat{A} + \hat{B}K\right)\Sigma \left(\hat{A} + \hat{B}K\right)^{\top}.
	\end{equation} 
	The problem formulation \eqref{prob:indirect}-\eqref{equ:indirect_sigma} is termed {\em indirect certainty-equivalence}  LQR design \cite{dorfler22on}. 
	

	\subsection{Direct certainty-equivalence LQR with covariance parameterization}
	Instead of estimating a dynamical model \eqref{equ:sysid}, the{ \em direct data-driven} LQR design aims to find $K^*$ directly from a batch of PE data~\cite{de2019formulas,dorfler2021certainty,zhao2024data}. For this purpose, our previous work \cite{zhao2024data} proposes a policy parameterization based on the sample covariance of input-state data
	\begin{equation}\label{def:cov}
	\Phi :  = \frac{1}{t}D_0D_0^{\top}=
	\begin{bmatrix}
	U_0D_0^{\top}/t  \\
	 \hline X_0D_0^{\top}/t 
	\end{bmatrix}
	=\begin{bmatrix}
	\overline{U}_0 \\
 \hline	\overline{X}_0
	\end{bmatrix}.
	\end{equation} Under the PE rank condition (\ref{equ:rank}), the sample covariance $\Phi$ is positive definite, and there exists a \textit{unique} solution $V\in \mathbb{R}^{(n+m)\times n}$ to
	\begin{equation}\label{equ:newpara}
	\begin{bmatrix}
	K \\
	I_n
	\end{bmatrix}= \Phi V
	\end{equation}
	for any given $K$. We refer to \eqref{equ:newpara} as the \textit{covariance parameterization} of the policy. Note that the dimension of the parameterized policy $V$ does not scale with data length $t$.

	For brevity, define the data matrices
	$
	\overline{W}_0= W_0D_0^{\top}/t, \overline{X}_1= X_1D_0^{\top}/t,
	$ {whose dimensions do not depend on $t$.}
	Then, the closed-loop matrix can be written as
	\begin{equation}\label{equ:clos}
	A+BK=[B,A]\begin{bmatrix}
	K \\
	I_n
	\end{bmatrix}\overset{\eqref{equ:newpara}}{=}[B,A]\Phi V\overset{\eqref{equ:dynamics}}{=}(\overline{X}_1 - \overline{W}_0)V.
	\end{equation}
	Following the certainty-equivalence principle, we disregard the uncertainty  $\overline{W}_0$ in \eqref{equ:clos}. Substituting $A+BK$ with $\overline{X}_1V$ in (\ref{equ:transfer})-(\ref{equ:Sigma}) and together with \eqref{equ:newpara}, the LQR problem becomes  
	\begin{equation}\label{prob:equiV}
	\begin{aligned}
	&\mathop{\text {minimize}}\limits_{V}~J(V) :=\text{Tr}\left((Q+V^{\top}\overline{U}_0^{\top}R\overline{U}_0V)\Sigma\right)\\
	&\text{subject to}~ ~\overline{X}_0V= I_n
	\end{aligned}
	\end{equation}
	with the gain matrix $K = \overline{U}_0V$, where $\Sigma$ is the positive definite solution to the Lyapunov equation
	\begin{equation}\label{equ:directSigma}
	\Sigma = I_n + \overline{X}_1V\Sigma V^{\top}\overline{X}_1^{\top}.
	\end{equation}
	 The problem formulation (\ref{prob:equiV})-\eqref{equ:directSigma} is termed \textit{direct certainty-equivalence} LQR design, and its solution coincides with that of the indirect design \eqref{prob:indirect}-\eqref{equ:indirect_sigma} \cite{zhao2024data}.

	\section{Policy gradient adaptive control}\label{sec:PGAC}
 	 In this section, we formulate the adaptive control problem and propose policy gradient adaptive control with convergence and stability guarantees.
 	  
 	\subsection{Problem statement} 
 	We focus on the following adaptive control problem.
 	\begin{problem}\label{problem}
 		design a gradient-based method leveraging online
 		closed-loop data such that the control policy converges to the
 		optimal LQR gain, while ensuring the stability of the closed-loop system.
 	\end{problem}
  	
  	Problem \ref{problem} concerns both stability and optimality, which is in line with the definition of adaptive control by Zames \cite{zames1998adaptive}, i.e., improve over the best with prior information. Throughout the article, we make the following assumption.\footnote{From Section \ref{sec:PGAC}, we use $X_{0,t}, U_{0,t}, W_{0,t}, X_{1,t}$ to denote the data series of length $t$ in \eqref{equ:data}. We also add a subscript $t$ to other notations to highlight the time dependence. }.

  	\begin{assum}\label{ass:PE}
  		The control input is persistently exciting such that $\underline{\sigma}(\Phi_t) \geq \gamma_t $ for some $\gamma_t > 0$.
  	\end{assum}
  	
  	Assumption \ref{ass:PE} is a quantitative condition of persistency of excitation and implies the rank condition \eqref{equ:rank}. It can be satisfied by adding probing noise to the control input and holds under another PE definition~\cite[Definition 2]{coulson2022quantitative}, which can be verified by examining solely the inputs. Notice again that the PE assumption is universal in adaptive control~\cite{bian2016value, lewis2009reinforcement}.
  	
  	

  	By Zames \cite{zames1998adaptive}, adaptive control involves the acquisition of \textit{information} about the plant. To specify the information metric, let $$\delta_t:=\|W_{0,t}D_{0,t}^{\top}/t\|$$ denote the bound of the covariance between the process noise and input-state data. This quantity reflects not only the amount of noise but also its correlation with the input-state data, which can capture particular statistics (e.g., i.i.d. noise). Note that we do not assume statistics of noise for the seek of generality. Define the \textit{signal-to-noise ratio} (SNR) of data at time $t$ as
  	$$
  	\text{SNR}_t := \frac{\gamma_t}{\delta_t},
  	$$ 
  	which describes the ratio between the useful and useless information \cite{dorfler2021certainty,dorfler22on}. We adopt the SNR as an \textit{information metric}, as the model uncertainty scales inversely with $\text{SNR}_t$.
  	
  	
  	\begin{lemma}\label{lem:id_error}
  		Consider the identification problem \eqref{equ:sysid} with data $(X_{0,t}, U_{0,t}, X_{1,t})$ and denote the model estimate as $(\hat{A}_t, \hat{B}_t)$. Then, it holds that
  		$$
  		\left\|[\hat{B}_t, \hat{A}_t] - [B, A] \right\| \leq 
  		\frac{1}{ \text{SNR}_t}.
  		$$
  	\end{lemma}
  	 
  	This information metric satisfies \textit{Zames's first monotonicity principle} for adaptation and learning \cite{zames1998adaptive}, as $\text{SNR}_t$ is usually a monotone increasing function of time. 
  	
  	\begin{remark}\label{rem:SNR}
  	We discuss the information $\text{SNR}_t$ under different noise and excitation settings that are widely adopted in adaptive control~\cite{simchowitz2020naive,zhao2024data,lu2023almost}. For time-uniformly bounded noise ($W_{0,t}$ and $D_{0,t}$ may be correlated) $\delta_t \sim \mathcal{O}(1)$ and constant excitation $\gamma_t \sim \mathcal{O}(1)$, we have constant $\text{SNR}_t \sim \mathcal{O}(1)$. For i.i.d. Gaussian noise, it holds in the probabilistic sense that $\delta_t  \sim \mathcal{O}(1/\sqrt{t})$~\cite{lu2023almost}. In this case, we have $\text{SNR}_t \sim \mathcal{O}(\sqrt{t})$ for constant excitation $\gamma_t \sim \mathcal{O}(1)$ and $\text{SNR}_t \sim \mathcal{O}(t^{\frac{1}{4}})$ for diminishing excitation $r_t \sim \mathcal{O}(t^{-\frac{1}{4}})$. \qed
  	\end{remark}

  	 
%
%
 	 
 	\subsection{The policy gradient adaptive control framework}
 	 To solve Problem \ref{problem}, we propose a \textit{policy gradient adaptive control}  (PGAC) framework alternating between control and gradient-based policy update. Specifically, the control input is in the form of $u_t = K_t x_t + e_t$, where the state-feedback gain $K_t$ is the parameterized policy at time $t$, and $e_t$ is a probing noise ensuring the PE condition. We assume that the initial gain $K_{t_0}$ is stabilizing, which is common in adaptive control~\cite{mania_certainty_2019,lu2023almost,simchowitz2020naive,cohen2019learning,chekan2024fully,karafyllis2019adaptive}. Define $\mathcal{S}:=\{K\in \mathbb{R}^{m\times n}|\rho(A+BK)<1\}$ as the set of stabilizing gains. 
 	 \begin{assum}\label{assum:gain}
 	 	The initial gain is stabilizing, i.e., $K_{t_0} \in \mathcal{S}$.
 	 \end{assum}
 	 
 	 We let Assumption \ref{assum:gain} hold in the rest of the article.
 	 The policy $K_t$ is updated over time with gradient methods of the LQR cost. When $(A,B)$ are known, the policy iterates as
 	 \begin{equation}\label{equ:grad}
 	 	K_{t+1} = K_t - \eta \nabla C(K_t),
 	 \end{equation}
 	 where $\nabla C(K_t)$ is the exact policy gradient. Then, the closed-form expression of $\nabla C(K)$ is given below.
 	 \begin{lemma}[{\cite[Lemma 1]{fazel2018global}}]\label{lem:mb_pg}
 	 	For any $K\in\mathcal{S}$, the gradient of $C(K)$ is given by
 	 	\begin{equation}\label{equ:gradK}
 	 		\nabla C(K) =  2\left(\left(R+{B}^{\top} P {B}\right) K+{B}^{\top} P {A}\right) {\Sigma},
 	 	\end{equation} 
 	 	where ${\Sigma}$ satisfies \eqref{equ:Sigma}, and $P$ is the positive definite solution to the Lyapunov equation 
 	 	\begin{equation}\label{equ:Lyap}
 	 		P = Q + K^{\top}RK + ({A}+{B}K)^{\top}P ({A}+{B}K).
 	 	\end{equation}
 	 \end{lemma}
 	 
 	 By Lemma 1, given knowledge of $(A, B)$, the policy gradient of the LQR cost can be computed by solving two Lyapunov equations in \eqref{equ:Sigma} and \eqref{equ:Lyap}. While the cost $C(K)$ is non-convex, it satisfies favorable properties such as gradient dominance and local smoothness, leading to global linear convergence of the gradient descent \eqref{equ:grad}  \cite{fazel2018global,bu2019lqr}.
 	 
 	 	\begin{lemma}[Gradient dominance]\label{lem:pl}
 	 	For $K \in \mathcal{S}$, it holds that
 	 	$$
 	 	C(K)- C(K^*) \leq \mu \| \nabla C(K)\|_F^2,
 	 	$$
 	 	where
 	 	$
 	 	\mu := \|\Sigma^*\|/\underline{\sigma}(R)
 	 	$ is the gradient dominance constant, and $\Sigma^*$ is the closed-loop covariance matrix \eqref{equ:Sigma} with $K^*$.
 	 \end{lemma}
 	  
 	 \begin{lemma}[{Local smoothness}]\label{lem:smooth}
 	 	Define the sublevel set $\mathcal{S}(a):=\{K\in \mathcal{S}|C(K)\leq a\}$. For any $K,K' \in \mathcal{S}(a)$ satisfying $K+b(K'-K) \in \mathcal{S}(a), \forall b \in [0,1]$, it holds that
 	 	$$
 	 	C(K') \leq C(K) +  \langle \nabla C(K), K' - K \rangle + {l(a)}\|K'-K\|^2/2,
 	 	$$
 	 	where the smoothness constant $l(a)$ is a polynomial in $a$.
 	 \end{lemma}
 	 
 	 Even when all gains are stabilizing $K_t \in \mathcal{S}, \forall t \geq t_0$, due to switches of the policies, the stability of the closed-loop system is unclear. Moreover, since $(A, B)$ are unknown in the adaptive control setting, the policy gradient can only be approximated using online closed-loop data, and the corresponding convergence remains unclear. 
 	  
 	 Depending on how the policy gradient is approximated, we categorize PGAC as indirect, when the gradient is computed with an estimated model~\eqref{equ:sysid}, versus direct, when the gradient is directly computed from data with the covariance parameterization \eqref{equ:newpara}-\eqref{prob:equiV}. In the sequel, we propose   indirect and direct PGAC methods with convergence and stability guarantees.


%

 	\subsection{Indirect policy gradient adaptive control}\label{subsec:indirect}
 	
 	\begin{algorithm}[t]
 		\caption{Indirect Policy Gradient Adaptive Control}
 		\label{alg:indirect}
 		\begin{algorithmic}[1]
 			\Require Offline data $(X_{0,t_0}, U_{0,t_0}, X_{1,t_0})$, an initial stabilizing policy $K_{t_0}$, and a stepsize $\eta$.
 			\For{$t=t_0,t_0+1,\dots$}
 			\State Apply $u_t = K_tx_t + e_t$ and observe $x_{t+1}$.
 			\State Estimate the model via recursive least-squares \eqref{equ:id}.
 			\State Perform one-step policy gradient descent 
 			\begin{equation}\label{equ:mbpg}
 				K_{t+1} = K_t - \eta \nabla \hat{C}_{t+1}(K_t),
 			\end{equation}
 			where $\nabla \hat{C}_{t+1}(K_t)$ is the policy gradient \eqref{equ:gradK} with the estimated model $(\hat{A}_{t+1},\hat{B}_{t+1})$.
 			\EndFor	
 		\end{algorithmic}
 	\end{algorithm}
 	
 	Indirect PGAC uses the least-squares estimates $(\hat{A},\hat{B})$ in \eqref{equ:sysid} to approximately compute the policy gradient. The details are presented in Algorithm \ref{alg:indirect}. We require an initial stabilizing gain and a batch of offline PE data. In the online stage, the control input contains a probing noise $e_t$ to ensure Assumption \ref{ass:PE}. After observing the new state, we use \textit{recursive least-squares} to identify an updated model estimate.   Let $\phi_t = [u_t^{\top},x_t^{\top}]^{\top}$. Then, the recursive least-squares iteration is 
 		\begin{equation}\label{equ:id}
 		[\hat{B}_{t+1}, \hat{A}_{t+1}] = [\hat{B}_{t}, \hat{A}_{t}] +  (x_{t+1} - [\hat{B}_{t}, \hat{A}_{t}]\phi_t  ) \frac{\phi_t^{\top} \Phi_t^{-1} }{t+ \phi_t^{\top}\Phi_t^{-1}\phi_t},
 		\end{equation} 
 		which is an efficient rank-one update. Then, we perform one-step policy gradient descent with stepsize $\eta$, where the gradient is computed with the estimated model. 
 	
 	

    Lemma \ref{lem:id_error} will play an important role in quantifying the effects of replacing the exact policy gradient with the approximated policy gradient in \eqref{equ:mbpg}. Together with Lemmas \ref{lem:pl} and \ref{lem:smooth}, we can show the convergence of $K_t$ in Algorithm \ref{alg:indirect}.

 	Apart from the convergence, we need to show stability of the closed-loop system. Since $K_t$ is updated over time, we resort to stability analysis of switching or time-varying systems \cite{cohen2019learning}.

 	\begin{definition}[Strong stability] \label{def:ss}
 		A policy $K$ is $(\kappa, \alpha)$-strongly stable if there exist constants $\kappa \geq 1$, $0< \alpha \leq 1$, and matrices $H \succ 0$ and $L$, such that  $A +BK = HLH^{-1}$ with $\|L\|\leq 1-\alpha$, $\|H\|\|H^{-1}\| \leq \kappa$, and $\|K\|\leq \kappa$.
 	\end{definition}
 	
 	It is known that a policy $K \in \mathcal{S}$ if and only if $K$ is $(\kappa, \alpha)$-strongly stable for some $\kappa$ and $\alpha$~\cite{cohen2019learning}. Definition \ref{def:ss} is quantitative and will lead to explicit bounds on the state. Since a sequence of policies $\{K_t\}$ is applied to the system, we require a stronger notion of stability under switching policies~\cite{cohen2019learning}. 
 	\begin{definition}[Sequential stability] \label{def:sss}
 		A sequence of policies $\{K_t\}$ for the system \eqref{equ:sys} is sequentially stable if there exist constants $\kappa \geq 1$, $0< \alpha \leq 1$, and matrices $\{H_t\succ 0\}$ and $\{L_t\}, \forall t$, such that $A+BK_t = H_tL_tH_t^{-1}$, and
 		\begin{enumerate}[label = (\roman*)]
 			\item $\|L_t\|\leq 1 - \alpha$ and $\|K_t\|\leq \kappa$;
 			\item $\|H_t\|\leq \kappa$ and $\|H_t^{-1}\|\leq 1$;
 			\item $\|H_{t+1}^{-1}H_t\| \leq 1 + \alpha/2$.
 		\end{enumerate}
 	\end{definition}
 	
 	With sequential stability, we can show boundedness of the state \cite{cohen2019learning}. By Definition \ref{def:sss}, sequential stability holds if $\{K_t\}$ is  strongly stable \textit{uniformly in time}, and the change in two consecutive policies is sufficiently small. In PGAC, this can be achieved by carefully controlling the stepsize $\eta$. In the following theorem, we show that  Algorithm \ref{alg:indirect} achieves both stability and optimality.

	\begin{theorem}\label{thm:indirect}
	There exist constants $\nu_i>0,i\in\{1,2,\dots,5\}$ with $\nu_4 <1$ depending on $(A,B,Q,R,K_{t_0})$, such that, if $\text{SNR}_t \geq \nu_1,$ for all $t$ and $\eta \leq \min\{\nu_2, 2\mu\}$, then the sequence $\{K_t\}$ is sequentially stable, and  the state is bounded as
	\begin{equation}\label{equ:bound_state}
	\|x_t\| \leq   \nu_3\left(1-\frac{\nu_4}{2}\right)^{t-t_0}\|x_{t_0}\|+\frac{2  \nu_3}{\nu_4} \max_{t_0 \leq i<t}\|Be_i + w_i\|.
	\end{equation} 
	Moreover, the policy  $\{K_t\}$ converges in the sense that
	\begin{equation}\label{equ:conv}
		\begin{aligned}
		C(K_{t}) - C^* &\leq \left(1-\frac{\eta}{2\mu}\right)^{t-t_0} (C(K_{t_0}) - C^*) \\
		&+ \eta \nu_5\sum_{i = t_0}^{t}\left(1-\frac{\eta}{2\mu}\right)^{t-i} \frac{1}{\text{SNR}_i}.
		\end{aligned}
	\end{equation} 
	\end{theorem} 
	
	We make several remarks on Theorem \ref{thm:indirect}. First, under a sufficiently large SNR and a small stepsize, the system is sequentially stable under switches of $K_t$. The state is explicitly bounded by two terms, i.e., an exponential decrease in the initial state plus an upper bound on the probing and process noise. In comparison, stability of the adaptive dynamical programming (ADP) approach~\cite{lewis2009reinforcement, bian2016value} has not been certified and is unclear. While the one-shot-based adaptive control methods ensure stability, they require a dwell time to mitigate the burn-in effect of switching policies and a safety mechanism in case of divergence~\cite{ chekan2024fully, cohen2019learning, lu2023almost}.


	 Second, the optimality gap of the LQR gain is upper bounded by two terms signifying an exponential decrease in the initial optimality gap plus a bias scaling inversely with the SNR. This aligns with \textit{Zames's second monotonicity principle} \cite{zames1998adaptive}, i.e., the performance monotonously increases with the information metric. Our results capture the convergence rate of the policy for different SNR functions in Remark \ref{rem:SNR}. For example, the optimality gap decreases as $\mathcal{O}(1/\sqrt{t})$ for $\text{SNR}_t \sim \mathcal{O}(\sqrt{t})$, which corresponds to the case of Gaussian noise and a constant excitation level, and $\mathcal{O}(t^{-\frac{1}{4}})$ for the diminishing excitation case with $\text{SNR}_t \sim \mathcal{O}(t^{\frac{1}{4}})$. Notice that the corresponding known optimal convergence rates are $\mathcal{O}({1}/{t})$ and $\mathcal{O}({1}/\sqrt{t})$, respectively, obtained by one-shot-based adaptive control \cite{simchowitz2020naive}. The latter requires a more conservative SNR condition such that the perturbation of the Riccati equation has a quadratic local approximation \cite{mania_certainty_2019}. Nevertheless, we observe that Algorithm \ref{alg:indirect} achieves the optimal convergence rate in simulations.

%

 	\subsection{Direct policy gradient adaptive control}\label{subsec:direct}
    Instead of computing the policy gradient with an estimated model, direct PGAC updates the policy based on the sample covariance parameterization~\eqref{equ:newpara}-\eqref{prob:equiV}. The details are presented in Algorithm \ref{alg:deepo}, which is also known as data-enabled policy optimization (DeePO) \cite{zhao2024data}. In the online stage, we use sample covariance to parameterize the policy \eqref{equ:k2v} and perform one-step projected gradient descent on the parameterized policy \eqref{equ:pgdV}. Define $\mathcal{V}: = \{V\in \mathbb{R}^{(n+m)\times n} |\rho(\overline{X}_1V)<1 \}$ as the feasible set of the covariance-parameterized LQR \eqref{prob:equiV}. The gradient and projection expressions are as below \cite{zhao2024data}.
     
    \begin{lemma} \label{lem:deepo_pg}
    	For $V \in \mathcal{V}$, the gradient of $J(V)$ is  
    	\begin{equation}\label{equ:pg}
    		\nabla J(V) = 2 \left(\overline{U}_{0}^{\top}R\overline{U}_{0}+\overline{X}_1^{\top}P\overline{X}_{1}\right)V \Sigma,
    	\end{equation}
    	where $\Sigma$ is the solution to \eqref{equ:directSigma}, and $P$ is the positive definite solution to the Lyapunov equation 
    	\begin{equation}\label{equ:Lya_V}
    		P = Q + V^{\top}\overline{U}_{0}^{\top}R\overline{U}_{0}V + V^{\top}\overline{X}_{1}^{\top}P\overline{X}_{1}V.
    	\end{equation} 
    	Further, we define the projection $\Pi_{\overline{X}_{0}}: = I_{n+m}-\overline{X}_{0}^{\dagger}\overline{X}_{0}$ onto the  constraint $\overline{X}_0V= I_n$ in \eqref{prob:equiV}.
    \end{lemma}
     
	By Lemma \ref{lem:deepo_pg}, given data $(X_0, U_0, X_1)$, the policy gradient of the LQR cost can be computed by solving two Lyapunov equations in \eqref{equ:directSigma} and \eqref{equ:Lya_V}. The projection is adopted to ensure the subspace constraint in \eqref{prob:equiV}. The stepsize $\eta_t$ is a variable to be specified later. As in indirect PGAC, Algorithm \ref{alg:indirect} also has an efficient recursive implementation \cite{zhao2024data}.  
    
    As a novel result beyond \cite{zhao2024data}, we show that the direct policy update \eqref{equ:k2v}-\eqref{equ:KV'} via projected gradient and the indirect counterpart \eqref{equ:mbpg} are equivalent up to a scaling matrix. 
    \begin{lemma}\label{lem:equi}
    	The policy update \eqref{equ:k2v}-\eqref{equ:KV'} is equivalent to
    	\begin{equation}\label{equ:relation}
    	K_{t+1} = K_t - \eta_t  M_{t+1} \nabla \hat{C}_{t+1}(K_t),
    	\end{equation} 
    	where $M_t:= \overline{U}_{0,t} \Pi_{\overline{X}_{0,t}} \overline{U}_{0,t}^{\top}$ is a positive definite matrix depending only on data. Moreover, the minimal singular value of $M_t$ satisfies $\underline{\sigma}(M_t) \geq {\gamma_t^2}$. 
    \end{lemma}
%

    Lemma \ref{lem:equi} enables us to leverage the previous results of indirect PGAC for the analysis of Algorithm \ref{alg:deepo}.

    \begin{algorithm}[t]
    	\caption{Direct Policy Gradient Adaptive Control}
    	\label{alg:deepo}
    	\begin{algorithmic}[1]
    		\Require Offline data $(X_{0,t_0}, U_{0,t_0}, X_{1,t_0})$, an initial stabilizing policy $K_{t_0}$, and stepsizes $\eta_t$.
    		\For{$t=t_0,t_0+1,\dots$}
    		\State Apply $u_t = K_tx_t + e_t$ and observe $x_{t+1}$.
    		\State Given $K_{t}$, solve $V_{t+1}$ via 
    		\begin{equation}\label{equ:k2v}
    			V_{t+1} =\Phi_{t+1}^{-1} \begin{bmatrix}
    				K_{t} \\
    				I_n
    			\end{bmatrix}.
    		\end{equation} 
    		\State Perform one-step projected gradient descent
    		\begin{equation}\label{equ:pgdV}
    			V_{t+1}' = V_{t+1} - \eta_{t+1}  \Pi_{\overline{X}_{0,t+1}} \nabla J_{t+1}(V_{t+1}),
    		\end{equation}
    		where the gradient $\nabla J_{t+1}(V_{t+1})$ is given by Lemma \ref{lem:deepo_pg}.
    		\State Update the control gain by 
    		\begin{equation}\label{equ:KV'}
    			K_{t+1} = \overline{U}_{0,t+1}V_{t+1}'.
    		\end{equation} 
    		\EndFor	
    	\end{algorithmic}
    \end{algorithm}
    
    \begin{theorem}\label{thm:direct}
    	There exist constants $\nu_i>0,i\in\{1,2,\dots,6\}$ with $\nu_4<1$ depending on $(A,B,Q,R,K_{t_0})$, such that, if $\text{SNR}_t \geq \max\{\nu_1, \nu_2\|M_{t}\|/\underline{\sigma}(M_{t}) \} $ and $\eta_t \leq  \|M_{t}\|^{-1} \min \{\nu_5, 2\mu\}$ for all $t$, then $\{K_t\}$ is sequentially stable, and \eqref{equ:bound_state} holds.
    	Moreover, the policy  converges as
    	   	\begin{align*}
    		C(K_{t}) - C^* &\leq \prod_{i=t_0+1}^{t}\left(1-\frac{\eta_i\underline{\sigma}(M_{i})}{2\mu}\right)(C(K_{t_0}) - C^*) \\
    		&+ \nu_6\sum_{i = t_0}^{t}\prod_{s=i}^{t-1}\left(1-\frac{\eta_s\underline{\sigma}(M_s)}{2\mu}\right) \frac{1}{\text{SNR}_i}.
    		\end{align*} 
    \end{theorem} 
    
  Theorem \ref{thm:direct} has three differences with Theorem \ref{thm:indirect} due to the relation between indirect and direct policy update in Lemma \ref{lem:equi}. First, the allowed SNR depends on the condition number of the scaling matrix $M_t$. Second, the stepsize $\eta_t$ scales with the inverse of $\|M_t\|$. Third, the convergence rate of the policy depends on the minimal singular value of $M_t$, which by Lemma \ref{lem:equi} scales with the excitation level $\gamma_t$. For the case of constant excitation, $\underline{\sigma}(M_t)$ is lower bounded by a constant, and the convergence of the policy reduces to that of indirect PGAC \eqref{equ:conv}. The diminishing excitation case does not lead to an explicit rate as $\underline{\sigma}(M_t)$ is time-varying. Note that in adaptive control it is essential to maintain sufficient excitation to prevent the bursting phenomenon \cite{anderson2005failures}. We refer to more comparisons between direct and indirect PGAC to simulations in Section \ref{sec:simu}.

  We discuss the special case of uniformly bounded noise and constant excitation considered in our previous work \cite{zhao2024data}. In this case, the minimal singular value of $M_t$ is lower bounded by a constant, and both $\|M_t\|$ and the condition number of $M_t$ converge asymptotically. Thus, the effect of $M_t$ can be neglected, and we can take a constant stepsize $\eta$. Then, the condition in Theorem \ref{thm:direct} reduces to that of  \cite[Theorem 2]{zhao2024data}. Theorem \ref{thm:direct} reveals a linear convergence of the policy with respect to the initial optimality gap, which improves over the sublinear rate  previously shown in \cite[Theorem 2]{zhao2024data}.

    \section{PGAC with the natural gradient and Gauss-Newton methods}\label{sec:grad}
    In the previous section, we have developed indirect and direct PGAC methods, where the vanilla  gradient descent is used to update the policy. In this section, we demonstrate the merits of two variants of gradient descent: natural gradient and Gauss-Newton methods, for the policy update of PGAC. In particular, the natural gradient descent bridges indirect and direct PGAC, and the Gauss-Newton method of the indirect PGAC leads to an adaptive version of Hewer's algorithm \eqref{equ:hewer}.
    
 	\subsection{Bridging indirect and direct PGAC via natural gradient}
 	
 	
 	Let us vectorize $K$ as $\theta = \text{vec}(K^{\top})$. Then, the \textit{natural policy gradient} update of the LQR cost $C(\theta)$ is given by\footnote{With a slight abuse of notation, we use $C(\theta)$ and $C(K)$ interchangeably.}
 	\begin{equation}\label{equ:ng}
 		\theta' = \theta - \eta F_{\theta}^{-1}\nabla C(\theta),
 	\end{equation} 
 	where $\nabla C(\theta)$ is the vanilla gradient, and  $F_{\theta}\succ 0$ is the Fisher information matrix (FIM) that captures the curvature of the parameter space of $\theta$, which biases the gradient descent in the direction of large uncertainties. The corresponding update for $K$ takes the form \cite{fazel2018global}
 	\begin{equation}\label{equ:ngK}
 		\begin{aligned}
 			K' &= K - \eta \nabla C(K)\Sigma^{-1} \\
 			& = K- 2\eta(\left(R+{B}^{\top} P {B}\right) K-{B}^{\top} P {A}),
 		\end{aligned}
 	\end{equation} 
 	where ${\Sigma}$ satisfies \eqref{equ:Sigma}, $P$ satisfies \eqref{equ:Lyap}, and the last equality follows from  Lemma \ref{lem:mb_pg}. We briefly highlight the merits of natural gradient.
 	

 	\begin{remark}\label{remark:na}
 		Compared with the vanilla gradient descent, the natural gradient descent \eqref{equ:ngK} is computationally more efficient, as it can be computed with only one Lyapunov equation \eqref{equ:Lyap} instead of two for the vanilla gradient  in Lemma \ref{lem:mb_pg}. Moreover, the stepsize for natural gradient \eqref{equ:ngK} can be chosen more aggressively than that of the vanilla gradient, leading to an improved  convergence rate \cite{fazel2018global}.  \qed
 	\end{remark} 	
 	
 	Motivated by \eqref{equ:ngK}, for indirect PGAC with natural gradient, we substitute the policy update \eqref{equ:mbpg} in Algorithm \ref{alg:indirect} with
 	\begin{equation}\label{equ:ngK_indirect}
 			K_{t+1} = K_t - 2\eta((R+\hat{B}_{t+1}^{\top} \hat{P}_{t+1} \hat{B}_{t+1}) K_t-\hat{B}_{t+1}^{\top} \hat{P}_{t+1} \hat{A}_{t+1}),
 	\end{equation}
 	where $\hat{P}_{t+1}$ is the positive definite solution to Lyapunov equation \eqref{equ:Lyap} with the estimated model $(\hat{A}_{t+1}, \hat{B}_{t+1})$.
 
 	For direct PGAC with natural gradient, we notice that the covariance parameterization \eqref{equ:newpara} defines a smooth bijection between the spaces of $K$ and $V$. By leveraging the fact that FIM captures the geometry of the paramterization space, we show that indirect and direct PGAC are bridged by
 	 adopting the natural gradient descent for the policy update. 
 	
 	 
 	 
 	\begin{lemma}\label{lem:bridge}
 	 For Algorithm \ref{alg:deepo}, 
 	 the policy update \eqref{equ:k2v}-\eqref{equ:KV'} with \eqref{equ:pgdV} replaced by natural gradient descent is equivalent to \eqref{equ:ngK_indirect}.
 	\end{lemma}
 	
 	We provide convergence and stability guarantees for PGAC with natural gradient.
 	\begin{theorem}\label{thm:na}
 		With different values of $\nu_i, i \in \{1,2,\dots,5\}$, the results in Theorem \ref{thm:indirect} also apply to PGAC with the natural gradient descent \eqref{equ:ngK_indirect}.
 	\end{theorem}
 	
    The theoretical guarantees in Theorem \ref{thm:na} are similar to those of PGAC with the vanilla gradient. Nevertheless, we note that performing natural gradient descent is computationally more efficient and leads to faster convergence rate; see Remark \ref{remark:na}. 
    

 	\subsection{Indirect PGAC with Gauss-Newton method and its equivalence to an adaptive Hewer's algorithm}
 		\begin{algorithm}[t]
 		\caption{Adaptive Hewer's Algorithm}
 		\label{alg:Hewer}
 		\begin{algorithmic}[1]
 			\Require Offline data $(X_{0,t_0}, U_{0,t_0}, X_{1,t_0})$ and an initial stabilizing policy $K_{t_0}$.
 			\For{$t=t_0,t_0+1,\dots$}
 			\State Apply $u_t = K_tx_t + e_t$ and observe $x_{t+1}$.
 			\State Estimate the model   via recursive least-squares \eqref{equ:id}.
 			\State Policy evaluation
 		 	$
 		 		\hat{P}_{t+1} = Q + K_t^{\top}RK_t 
 		 		+(\hat{A}_{t+1} + \hat{B}_{t+1}K_t)^{\top}\hat{P}_{t+1}(\hat{A}_{t+1} + \hat{B}_{t+1}K_t).
 		 	$ 
 		 	\State Policy improvement
 		 	$$
 		 	K_{t+1} = (R + \hat{B}_{t+1}\hat{P}_{t+1}\hat{B}_{t+1})^{-1}\hat{B}_{t+1}\hat{P}_{t+1}\hat{A}_{t+1}.
 		 	$$
 			\EndFor	
 		\end{algorithmic}
 	\end{algorithm}

 	The \textit{Gauss-Newton method} is a quasi-Newton method under a particular Riemannian metric \cite{bu2019lqr}. It iterates as
 	\begin{equation}\label{equ:pure_GN}
 		\begin{aligned}
 		 K' &= K - \eta \left(R+B^{\top}PB\right)^{-1}\nabla C(K)\Sigma^{-1}\\
 		&= K - 2\eta\left(R+B^{\top}PB\right)^{-1} ((R+{B}^{\top} P {B}) K-{B}^{\top} P {A} )
 		\end{aligned}
 	\end{equation}
 	with the constant stepsize $0<\eta <1$. The stepsize of the Gauss-Newton method can be chosen more aggressively than the vanilla and natural gradient, leading to a faster convergence rate. With a special choice of stepsize $\eta = 1/2$, the Gauss-Newton update \eqref{equ:pure_GN} is equivalent to
 	$$
 	K' = \left(R+B^{\top}PB\right)^{-1}B^{\top}PA.
 	$$
 	This coincides with the Hewer's algorithm \eqref{equ:hewer} and enjoys local quadratic convergence. However, the Gauss-Newton update is computationally less efficient than natural gradient, as computing \eqref{equ:pure_GN} involves the inverse of a matrix.  
 	
 	Inspired by \eqref{equ:pure_GN}, for indirect PGAC we substitute the policy update \eqref{equ:mbpg} in Algorithm \ref{alg:indirect} with 
 	\begin{equation}\label{equ:GNK_indirect}
 		\begin{aligned}
 			K_{t+1} &= K_t - 2\eta(R+\hat{B}_{t+1}^{\top}\hat{P}_{t+1}\hat{B}_{t+1})^{-1} \\
 			&\times((R+\hat{B}_{t+1}^{\top} \hat{P}_{t+1} \hat{B}_{t+1}) K_t-\hat{B}_{t+1}^{\top} \hat{P}_{t+1} \hat{A}_{t+1}).
 		\end{aligned} 
 	\end{equation}
 	Our theoretical guarantees of indirect PGAC with the Gauss-Newton update \eqref{equ:GNK_indirect} are as below.
 		\begin{theorem}\label{thm:GN}
 		With different values of $\nu_i, i \in \{1,2,\dots,5\}$ and $\eta < 2\|\Sigma^*\|$, the results in Theorem \ref{thm:indirect} apply to indirect PGAC with the Gauss-Newton policy update \eqref{equ:GNK_indirect}, with the convergence certificate
 		\begin{equation}\label{equ:conv_gnn}
 			\begin{aligned} 
 				C(K_{t}) - C^* &\leq \left(1-\frac{\eta}{2\|\Sigma^*\|}\right)^{t-t_0} (C(K_{t_0}) - C^*) \\
 				&+ \eta \nu_5\sum_{i = t_0}^{t}\left(1-\frac{\eta}{2\|\Sigma^*\|}\right)^{t-i} \frac{1}{\text{SNR}_i}.
 			\end{aligned}
 		\end{equation}  
 	\end{theorem}  
 	
 	Since the stepsize can be chosen aggressively, it usually leads to faster convergence than PGAC with the vanilla and natural gradient. With the stepsize $\eta = 1/2$, our indirect PGAC with Gauss-Newton update is equivalent to an adaptive version of the Hewer's algorithm (also referred to as online identification-based policy iteration in  \cite{song2025robustness}); see Algorithm \ref{alg:Hewer}. However, Algorithm \ref{alg:Hewer} generally does not have sequential stability guarantees, as the choice of $\eta = 1/2$ may not ensure a sufficiently small policy change. An exception is when the initial policy is sufficiently close to the optimal LQR gain.  
  \begin{theorem}\label{thm:indirect_gn}
 	 There exist constants $\nu_i>0,i\in\{1,2,\dots,5\}$ with $\nu_4<1$ depending on $(A,B,Q,R,K_0)$, such that, if $\text{SNR}_t \geq \nu_1, \forall t$ and $C(K_{t_0})-C^*\leq \nu_2$, then $\{K_t\}$ is sequentially  stable, and \eqref{equ:bound_state} holds. Moreover, the policy  $\{K_t\}$ converges in the sense that
 		\begin{equation}\label{equ:conv_gn_qd}
		\begin{aligned}
			C(K_{t+1}) - C^* &\leq  \frac{1}{\nu_2} (C(K_{t}) - C^*)^2 + \frac{ \nu_5}{\text{SNR}_{t+1}}.
		\end{aligned}
	\end{equation} 
 	\end{theorem} 
 	
 	
 	In Theorem \ref{thm:indirect_gn}, the condition on the initial policy $K_{t_0}$ can be replaced with a condition on the SNR, given that $K_{t_0}$ is obtained as the solution of the indirect certainty-equivalence LQR \eqref{prob:indirect} with offline data. This aligns with the convergence and stability results in \cite{song2025robustness}that hold under sufficiently small identification error. The local quadratic convergence in \eqref{equ:conv_gn_qd} inherits from properties of the quasi-Newton method~\cite{bu2019lqr}. 
 	Note that the dependence of the convergence rate on the SNR is the same as that of Theorem \ref{thm:indirect}.

 	\section{Boosting the performance of PGAC via regularization}\label{sec:reg} 
 	So far, we have developed PGAC with variants of gradient descent of the certainty-equivalence LQR formulations \eqref{prob:indirect} and \eqref{prob:equiV}. However, the certainty-equivalence approach neglects the effect of noise in the data, leading to uncertainties in the closed-loop covariance matrix and the cost function. This section recalls a regularization method for PGAC to compensate the uncertainty~\cite{zhao2025regularization}. As a main contribution, we show that the theoretical guarantees of PGAC are preserved under proper choices of the regularization coefficient.

\subsection{Indirect PGAC with regularization}
For the indirect certainty-equivalence LQR \eqref{prob:indirect}-\eqref{equ:indirect_sigma}, the closed-loop covariance matrix $\Sigma$ satisfies the Lyapunov equation \eqref{equ:indirect_sigma}.
However, given $(A, B)$, the Lyapunov equation that should be met is \eqref{equ:Sigma}. For brevity, define 
\[
\Xi := \begin{bmatrix}
	K \\ I_n
\end{bmatrix}\Sigma \begin{bmatrix}
	K \\ I_n
\end{bmatrix}^{\top}.
\]
Then, the difference between the right-hand sides of \eqref{equ:indirect_sigma} and \eqref{equ:Sigma} is 
\begin{align*}
	 \overline{W}_0 \Phi^{-1} \Xi [B~ A]^{\top} + [B~ A] \Xi \Phi^{-1} \overline{W}_0^{\top} + \overline{W}_0 \Phi^{-1} \Xi \Phi^{-1} \overline{W}_0^{\top}.
\end{align*}
For well-behaved noise statistics, $\|\overline{W}_0\|$ usually vanishes quickly with time \cite[Lemma 1]{zhao2025regularization}, and the first two terms dominate the difference. To reduce the difference, we introduce the regularizer $\text{Tr}(\Phi^{-1}\Xi)$ to the certainty-equivalence cost \eqref{prob:indirect}, leading to the regularized cost
\begin{equation}\label{equ:reg_cost}
	\hat{C}(K;\lambda) := \text{Tr}\left( \left(\begin{bmatrix}
		R & 0 \\
		0 & Q
	\end{bmatrix} + \lambda\Phi^{-1} \right) \begin{bmatrix}
		K \\
		I_n
	\end{bmatrix} \Sigma \begin{bmatrix}
		K \\
		I_n
	\end{bmatrix}^{\top}  \right),
\end{equation}
where $\Sigma$ satisfies \eqref{equ:indirect_sigma}, and $\lambda>0$ is the regularization coefficient. The regularizer is a correction to the penalty matrices of the LQR problem. In line with \cite{chiuso2023harnessing}, the regularization  also compensates the uncertainty in the cost function and promotes exploitation \cite{zhao2025regularization}. We demonstrate these merits via simulations in Section \ref{subsec:simu_reg}. 

Next, we derive the gradient expression of the regularized cost \eqref{equ:reg_cost}. Consider the following partition of $\Phi^{-1}$
$$
\Phi^{-1} = \begin{bmatrix}
	\left(\Phi^{-1}\right)_{uu} & \left(\Phi^{-1}\right)_{ux} \\
	\left(\Phi^{-1}\right)_{xu} & \left(\Phi^{-1}\right)_{xx}
\end{bmatrix},
$$
where $(\Phi^{-1})_{uu}\in \mathbb{R}^{m\times m}$ and $(\Phi^{-1})_{xx}\in \mathbb{R}^{n\times n}$. Define $Q_{\lambda} = Q + \lambda (\Phi^{-1})_{xx}$ and $R_{\lambda} = R + \lambda (\Phi^{-1})_{uu}$. 


\begin{lemma}\label{lem:grad_reg}
	For $K\in \{K|\rho(\hat{A}+\hat{B}K)<1\}$, the gradient of $\hat{C}(K;\lambda)$ is given by
	$$
	\nabla \hat{C}(K;\lambda) = 2\left(R_{\lambda}K + \hat{B}^{\top}P(\hat{A}+\hat{B}K) +\lambda \left(\Phi^{-1}\right)_{ux} \right)\Sigma,
	$$
	where $\Sigma$  and $P$ are the solutions to \eqref{equ:indirect_sigma} and
	$
	P = Q_{\lambda}  + K^{\top}R_{\lambda}K + \lambda K^{\top}(\Phi^{-1})_{ux} + \lambda(\Phi^{-1})_{xu}K + (\hat{A}+\hat{B}K)^{\top}P (\hat{A}+\hat{B}K)
	$, respectively.
\end{lemma} 

The proof of Lemma \ref{lem:grad_reg}  follows the same vein as that of \cite[Lemma 1]{fazel2018global} and is omitted due to space limitation. 

To apply regularization to the indirect PGAC, we replace the gradient in Algorithm \ref{alg:indirect} with that of the regularized cost in Lemma \ref{lem:grad_reg}. To ensure convergence of ${K_t}$ to $K^*$, the effect of the regularizer must vanish over time. Since $\text{Tr}(\Phi_t^{-1}\Xi)$ scales with $\|\Phi_t^{-1}\| \leq \mathcal{O}(1/\gamma_t)$, setting $\lambda_t \leq \mathcal{O}(\delta_t)$ ensures that the regularization decays as $\mathcal{O}(1/\text{SNR}_t)$. We show this decay rate is sufficient for convergence and stability of the regularized indirect PGAC.

\begin{theorem}\label{thm:indirect_reg}
Consider the regularized indirect PGAC algorithm, where we replace the gradient $\nabla\hat{C}_{t+1}(K_t)$ in \eqref{equ:mbpg} with $\nabla\hat{C}_{t+1}(K_t;\lambda_{t+1})$. With different values of $\nu_i, i \in \{1,2,\dots,5\}$ and  the condition $\lambda_t \leq \nu_6 \delta_t$ for $\nu_6>0$, the results in Theorem \ref{thm:indirect} apply, i.e., the boundedness of the state \eqref{equ:bound_state} and the convergence certificate of $C(K_t)$ \eqref{equ:conv} hold.
\end{theorem} 

Under the condition $\lambda_t \leq \mathcal{O}(\delta_t)$, the theoretical guarantees can be extended to PGAC with natural gradient and Gauss-Newton method, whose formal proofs are omitted.




 	\subsection{Direct PGAC with regularization} 
    Leveraging the covariance parameterization in \eqref{equ:newpara}, the regularizer in \eqref{equ:reg_cost} can be reformulated as $\text{Tr}(V\Sigma V^{\top}\Phi)$. Then, the direct LQR cost with regularization is 
 	\begin{equation}\label{prob:regu}
	 			  J(V;\lambda):= J(V) + \lambda\text{Tr}(V\Sigma V^{\top}\Phi).
	\end{equation}
	We derive the gradient expression of $J(V;\lambda)$, the proof of which follows the same vein as that of \cite[Lemma 2]{zhao2024data} and is omitted due to space limitation.
	\begin{lemma}\label{lem:gradient}
		For $V\in \mathcal{V}$, the gradient of $J(V;\lambda)$ is
		$$
		\nabla J(V;\lambda) = 2 \left(\lambda \Phi + \overline{U}_0^{\top}R\overline{U}_0+\overline{X}_1^{\top}P\overline{X}_1\right)V \Sigma,
		$$
		where $\Sigma$ and $P$ are the solution to \eqref{equ:directSigma} and  
		$
		P = Q  + V^{\top}(\lambda \Phi + \overline{U}_0^{\top}R\overline{U}_0)V + V^{\top}\overline{X}_1^{\top}P\overline{X}_1V
		$, respectively.
	\end{lemma}
%
 
	As in regularized indirect PGAC, the regularization coefficient should scale as $\lambda_t \leq \mathcal{O}(\delta_t)$. Then, by leveraging the relation between indirect and direct PGAC in Lemma \ref{lem:equi}, we have the convergence and stability results. 
 
	 \begin{theorem}\label{thm:direct_reg}
	 	Consider the regularized direct PGAC algorithm, where we replace the gradient $ \nabla J_{t+1}(V_{t+1})$ in \eqref{equ:pgdV} with $ \nabla J_{t+1}(V_{t+1};\lambda_{t+1})$ in Lemma \ref{lem:gradient}. With different values of $\nu_i, i \in \{1,2,\dots,6\}$ and the condition $\lambda_t \leq \nu_7 \delta_t$ for $\nu_7>0$, the boundedness of the state and the convergence certificate of $C(K_t)$  in Theorem \ref{thm:direct} apply.
	 \end{theorem}

	 We discuss the selection of the regularization coefficient $\lambda_t = \lambda_0 \delta_t$ under special cases of Remark \ref{rem:SNR}. For the bounded noise case, a constant $\lambda_t$ is sufficient. For the i.i.d. Gaussian noise case, we let $\lambda_t \sim \mathcal{O}(1/\sqrt{t})$. Our results align with \cite{zeng2024noise}, which shows that a constant $\lambda_t$ over time is not favorable and may lead to a trivial solution in direct parameterization \cite{de2019formulas}. Note that we still need to tune the coefficient $\lambda_0$ for optimal performance and stability.

	\section{Numerical case studies}\label{sec:simu}
	This section uses simulations to demonstrate the effectiveness of PGAC in terms of convergence, stability, and computational efficiency. Our simulations are based on the benchmark LQR problem with the model \cite[Section 6]{dean2020sample}
	\begin{equation}\label{equ:benchmark}
		\begin{aligned}
			&A = \begin{bmatrix}
				1.01 &  0.01  &  0\\
				0.01  & 1.01 &  0.01\\
				0  &  0.01  &  1.01
			\end{bmatrix}, ~B = I_3,
		\end{aligned}
	\end{equation}
	which corresponds to a discrete-time marginally unstable Laplacian system, and $Q= I_3$ and $R = 10^{-3} \times I_3$. Note that the direct PGAC method (also known as DeePO \cite{zhao2024data}) has been validated in nonlinear power systems \cite{zhao2024direct} and real-world experiments of an autonomous bicycle \cite{persson2025adaptive}. Our code is run in Matlab and provided in \url{https://github.com/Feiran-Zhao-eth/policy-gradient-adaptive-control}.
	
	\subsection{Convergence and closed-loop stability of PGAC}\label{subsec:conv_simu}
	
	\begin{figure}[t]
		\centerline{\includegraphics[width=70mm]{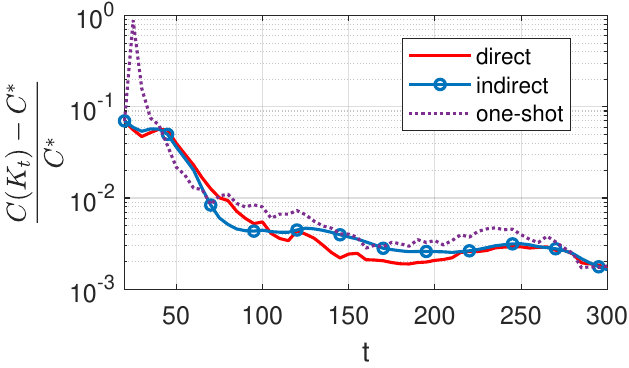}}
		\caption{Convergence of one-shot-based and PGAC methods.}
		\label{pic:convergence}
	\end{figure}
	
	We compare the convergence and stability of the one-shot-based adaptive control, direct and indirect PGAC with the vanilla gradient, natural gradient, and Gauss-Newton methods. We set the number of offline data to $t_0=20$ and let $x_0 = 0$, $u_t\sim \mathcal{N}(0,I_3)$ for $t<t_0$. For $t \geq t_0$, we set $u_t =K_tx_t + e_t$ with $e_t \sim \mathcal{N}(0,I_3)$ to ensure persistency of excitation. The noise is drawn from $w_t\sim  \mathcal{N}(0, I_3)$, which corresponds to a poor $\text{SNR}_t \in [-5,0]$ dB. We set the stepsize to $\eta = 0.02$ for  indirect PGAC with the vanilla gradient, $\eta = 0.2$ for natural gradient, $\eta = 0.5$ for the Gauss-Newton method (Algorithm \ref{alg:Hewer}), and $\eta_t = 0.2/\|M_t\|$ for direct PGAC. The one-shot-based method obtains $K_t$ from the Riccati equation with an estimated model \eqref{equ:id} every time step. For all the methods, we set $K_{t_0}$ as the certainty-equivalence LQR \eqref{prob:indirect} with offline data.
	
	Fig. \ref{pic:convergence} shows the optimality gap of the policy sequence for the one-shot-based method and direct and indirect PGAC with the vanilla gradient under the same randomness realization. For clarity of presentation, we omit the curves of indirect PGAC with natural gradient and Gauss-Newton, which lie between those of the one-shot-based method and indirect PGAC with the vanilla gradient.  We make several key observations. First, all the methods improve the performance over the initial policy using online closed-loop data, and the optimality gap converges asymptotically at the rate $\mathcal{O}(1/t)$, which implies that our theoretically certified rate $\mathcal{O}(1/\sqrt{t})$ is conservative. Second, the one-shot-based method diverges at the early stage due to unexpected large noise. In contrast, PGAC methods exhibit more stable convergence. This is expected in presence of noise, as the gradient is more robust than the minimizer.

	\begin{figure}[t]
		\centerline{\includegraphics[width=70mm]{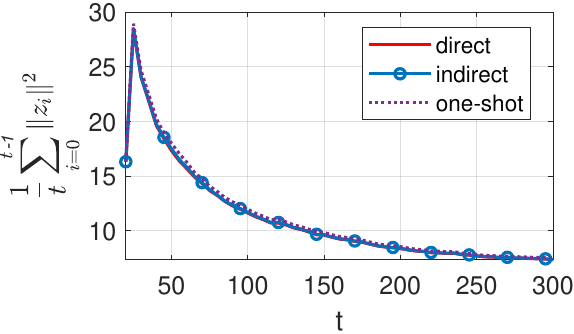}}
		\caption{Finite-horizon cost of one-shot-based and PGAC methods.}
		\label{pic:finite_cost}
	\end{figure}
	
 	Fig. \ref{pic:finite_cost} shows their average finite-horizon cost $\sum_{i=0}^{t-1} \|z_i\|^2/t$ as a function of time, where $z_t$ is the performance output in \eqref{equ:closedsys}. We observe that all methods exhibit highly similar average costs, which converge rapidly to a positive constant. This indicates that the closed-loop systems resulting from both PGAC and the one-shot-based method are stable. Furthermore, the convergence behavior aligns with our theoretical results, i.e., the state is bounded by the sum of an exponentially decaying term and a steady-state bias determined by the probing and process noise. We note that the divergence at the early stage is due to abrupt large process noise.
	
		\begin{table}[t]\label{table:time}
		\renewcommand{\arraystretch}{1.4}
		\begin{center}
			\begin{tabular}{|c|c|c|c|c|c|}
				\hline
				& direct & indirect & natural & GN & one-shot-based   \\ 
				\hline
				time (s)& 0.063  & 0.057	& 0.047 & 0.049 & 0.418   \\ 
				\hline
			\end{tabular}
		\end{center}
		\caption{Comparison of the computation time for running $500$ time steps.}
	\end{table}
	
	\begin{table}[t]\label{table}
		\renewcommand{\arraystretch}{1.4}
		\begin{center}
			\begin{tabular}{|c|c|c|c|}
				\hline
				& $\lambda_t = 0$ & $\lambda_t = 1/(10\sqrt{t-t_0})$ & one-shot-based   \\
				\hline
				\tabincell{c}{ direct }  & \tabincell{c}{ {$\mathcal{P} = 81\%$ } \\ {$\mathcal{M} = 0.0014$}} & \tabincell{c}{ {$\mathcal{P} = 98\%$ } \\ {$\mathcal{M} = 0.0011$}} & \tabincell{c}{ {$\mathcal{P} = 82\%$ } \\ {$\mathcal{M} = 0.0013$}} \\
				\hline
				\tabincell{c}{ indirect }  & \tabincell{c}{ {$\mathcal{P} = 83\%$ } \\ {$\mathcal{M} = 0.0014$}} & \tabincell{c}{ {$\mathcal{P} = 99\%$ } \\ {$\mathcal{M} = 0.0014$}} & \tabincell{c}{ {$\mathcal{P} = 82\%$ } \\ {$\mathcal{M} = 0.0013$}} \\
				\hline
			\end{tabular}
		\end{center}
		\caption{The impact of regularization for PGAC.}
	\end{table}
	
	\subsection{Computational efficiency of PGAC}
	We compare their efficiency in terms of the computation time. 
	For indirect PGAC, we apply recursive least
	squares (\ref{equ:id})  for identification. For direct PGAC, we apply the rank-one update in \cite[Section V-B]{zhao2024data}. For the one-shot-based method,  the certainty-equivalence LQR problem \eqref{prob:indirect} is solved using the \textit{dlqr} function in MATLAB. We perform $20$ independent trials and record the mean of the computation time for running $t=500$ time steps. The results are summarized in Table I. We observe that compared with the one-shot-based method, the PGAC approaches require significantly less computational time. This gap is due to that PGAC performs only one or two Lyapunov equation solves per iteration for gradient descent, whereas the one-shot-based method involves solving a Riccati equation at each step.  Among all the PGAC approaches, direct PGAC consumes highest computation time as the dimension of the optimization matrix $V\in \mathbb{R}^{(n+m)\times n}$ is higher than that of $K\in \mathbb{R}^{m\times n}$ in the indirect case. For indirect PGAC, the vanilla gradient method is slowest as it requires solving two Lyapunov equations per time step instead of one for natural gradient and Gauss-Newton methods. It is worth noting that while we demonstrate the results using a simple third-order system \eqref{equ:benchmark}, the computational gap between PGAC and the one-shot-based method grows substantially with increasing system dimension; we refer the reader to our previous work~\cite[Section V]{zhao2024data} for more comprehensive simulation results.

	\subsection{The impact of regularization}\label{subsec:simu_reg}
	
	We study the impact of regularization for both indirect and direct PGAC with the vanilla gradient. We set the stepsize to $\eta = 0.2$ for  indirect PGAC  and $\eta_t = 0.2/\|M_t\|$ for direct PGAC. Consider $100$ independent trials. We denote by $\mathcal{P}$ the percentage the algorithm converges without any instability issue and by $\mathcal{M}$ the median of the optimality gap ${(C(K_T) - C^*)}/{C^*}$ with $T=1000$ through all trails where convergence is observed. The results are reported in Table II, where the one-shot-based method without regularization is for comparison. With regularization, both indirect and direct PGAC significantly outperform the one-shot-based method in terms of algorithm stability, i.e., the robustness of policy updates against noise. Furthermore, direct PGAC achieves a smaller optimality gap compared to the one-shot-based method. This is because our regularizer well compensates the uncertainty induced by noise in both the closed-loop covariance matrix and the cost function~\cite{zhao2025regularization}.



	\section{Concluding remarks}\label{sec:conclu}
	In this article, we proposed policy gradient adaptive control (PGAC) for the LQR, which touches upon both indirect and direct approaches and enables the use of variant gradient methods and regularization. For all approaches, we showed convergence of the policy and stability of the closed-loop system. Simulations validated the theoretical results and illustrated the robustness and computational efficiency of PGAC.
	
	We believe that this article leads to fruitful future works.  It would be interesting to see if the probing signal can be designed using the sample covariance to efficiently reduce the uncertainty. As observed in the simulation, our theoretical convergence rate may be conservative, and a sharper analysis might be possible. Our sequential stability analysis can be used for other RL-based adaptive control methods, e.g., adaptive dynamical programming, where the closed-loop stability is largely open \cite{annaswamy2021historical}. Our PGAC framework can be extended to other settings (e.g., output feedback control), performance indices (e.g., $\mathcal{H}_{\infty}$ norm), noise assumptions (e.g., stochastic noise), and system classes (e.g., time-varying systems). 
	

%

	\appendices
	\section{Proofs in Section \ref{sec:PGAC}}
	 
	\subsection{Proof of Lemma \ref{lem:id_error}}
	It follows from Assumption \ref{ass:PE} that
	$
	\|[\hat{B}_t, \hat{A}_t] - [B, A]\| = 	\|W_{0,t} D_{0,t}^{\dagger}\| = \|\overline{W}_{0,t}\Phi_t^{-1}\| \leq {\delta_t}/{ \gamma_t}.
	$

	\subsection{Proof of Theorem \ref{thm:indirect}}	 
	We first provide some useful bounds.
	\begin{lemma}\label{lem:bounds}
		For $K \in \mathcal{S}$, it holds
		(\romannumeral1) 
		$
		\|\Sigma\|  \leq {C(K)}/{\underline{\sigma}(Q)},
		$
		(\romannumeral2) 
		$
		\|P\| \leq C(K),
		$
		and (\romannumeral3) 
		$
		\|K\|_F \leq ({C(K)}/{\underline{\sigma}(R)})^{{1}/{2}}.
		$
	\end{lemma}
	\begin{proof}
		The proof of (\romannumeral1) and (\romannumeral2) follows directly from the definition of $C(K)$. To show (\romannumeral3), we have
		$ C(K) = \text{Tr}((Q+K^{\top}RK)\Sigma)
		\geq \text{Tr}(K^{\top}RK\Sigma)
		\geq \text{Tr}(K^{\top}K) \underline{\sigma}(R) 
		= \|K\|_F^2 \underline{\sigma}(R).$
	\end{proof}	
	
	The proof frequently invokes the following perturbation analysis result for Lyapunov equations~\cite[Lemma 15]{zhao2024data}. 
	\begin{lemma}\label{lem:perturb}
		Let $A\in \mathbb{R}^{n \times n}$ be stable and $\Sigma(A)$ be the unique positive definite solution to $\Sigma(A) = I_n + A\Sigma(A) A^{\top}$. If 
		$
		\|A'-A\| \leq 1/({4\|\Sigma(A)\|(1+\|A\|)}),
		$
		then $A'$ is stable and
		$
		\|\Sigma(A')-\Sigma(A)\| \leq 4 \|\Sigma(A)\|^2(1+\|A\|)\|A'-A\| .
		$
	\end{lemma}

	Consider the policy gradient update with estimated $(\hat{A}, \hat{B})$
	\begin{equation}\label{equ:pg_appr}
		K' = K - \eta \nabla \hat{C}(K),
	\end{equation}
	and the update with the ground-truth $(A,B)$
	\begin{equation}\label{equ:pg_truee}
		K'' = K - \eta \nabla {C}(K).
	\end{equation}

	
	The convergence of the update  \eqref{equ:pg_truee} can be shown by leveraging  Lemmas \ref{lem:pl} and \ref{lem:smooth}. If $\eta \leq 1/l$, then
	\begin{equation}\label{equ:exact}
		C(K'') - C(K) \leq  -\frac{\eta}{2\mu} (C(K) - C^*), 
	\end{equation} 
	where $l$ is the smoothness parameter and $\mu$ is the gradient dominance constant.
	To show the convergence of \eqref{equ:pg_appr}, we first quantify the distance between the exact gradient and the approximated gradient. Let $p_1$ be a scalar function 
	$$
	p_1 = \frac{8C^3(K)}{\underline{\sigma}^2(Q)}\left(1+\frac{C(K)}{\underline{\sigma}(Q)}\right) \left(1+\sqrt{\frac{C(K)}{\underline{\sigma}(R)}}\right),
	$$
	and $p_2 = {C^2(K) }/({\underline{\sigma}(Q)p_1})$. For the ease of notation, we omit the subscript in $\gamma_t$ and $\delta_t$ when it is clear from the context.
	\begin{lemma}\label{lem:grad_diff}
		Let $K \in \mathcal{S}$. Then, there exists a polynomial  
		$p_3 = \text{poly}(C(K)/\underline{\sigma}(Q), \|A\|, \|B\|, \|R\|, 1/\underline{\sigma}(R))$ such that, if $\delta / \gamma \leq p_2$, then $
				 \| \nabla{C}(K) - \nabla\hat{C}(K) \| \leq  {p_3\delta}/{\gamma}.$
	\end{lemma}
	\begin{proof} 
		For brevity, let 
		\begin{equation}\label{equ:E}
			E = RK + B^{\top}P(A+BK).
		\end{equation} 
		We use $\hat{E}$, $\hat{\Sigma}$, and $\hat{P}$ to denote the quantities under the estimated model $(\hat{A}, \hat{B})$. Then, the gradient can be expressed as $\nabla C(K) = 2E\Sigma$, and the  difference in gradients is
		\begin{equation}\label{equ:diff}
			\nabla{C}(K) - \nabla\hat{C}(K) 
			=   + 2E(\hat{\Sigma} - \Sigma) + 2(\hat{E} - E) \hat{\Sigma}.
		\end{equation} 
		 
		We quantify the first term of \eqref{equ:diff}. By Lemma \ref{lem:id_error},
		it holds 
		\begin{equation}\label{equ:closeloop}
			\|(A+BK) - (\hat{A} + \hat{B}K)\|\leq (1+\|K\|){\delta}/{\gamma}.
		\end{equation}
		Then, we apply Lemma \ref{lem:perturb} to show 
		$
			\|\hat{\Sigma} - \Sigma\| \leq 4 \|\Sigma\|^2(1+\|A+BK\|)\|A+BK - (\hat{A} + \hat{B}K)\| 
			\leq  {p_1\delta}/{(2\gamma C(K))}.
		$ Together with the bound 
		$\operatorname{Tr}\left(E^{\top} E\right) \leq {\left\|R+B^{\top} P B\right\|\left(C(K)-C^*\right)}$ \cite[Lemma 11]{fazel2018global},
		the first term of \eqref{equ:diff} is bounded as
		$
		\|2E(\hat{\Sigma} - \Sigma)\|\leq {\left\|R+B^{\top} P B\right\|^{\frac{1}{2}} \left(C(K)-C^*\right)^{\frac{1}{2}}}p_1\delta/(2\gamma C(K)).
		$
		
		To bound the second term of \eqref{equ:diff}, we note that
		$
				E - \hat{E} = B^{\top}P((A+BK) - (\hat{A} + \hat{B}K)) 
				+ (B^{\top}P - \hat{B}^{\top}\hat{P})(\hat{A} + \hat{B}K).
		$
		The first term of the right-hand side of this equation can be bounded by using \eqref{equ:closeloop} and Lemma \ref{lem:bounds}. For the second term, we have 
		$\|\hat{A} + \hat{B}K\| \leq \|A+BK\| +  (1+\|K\|){p_2}$ and  $\|\hat{\Sigma}\|\leq \|\Sigma\| + {p_1p_2}/{(2C(K))}$. 
		Notice that  
		$
		B^{\top}P - \hat{B}^{\top}\hat{P} = \hat{B}^{\top}(P-\hat{P}) + (B-\hat{B})^{\top}{P}
		$.
		By the proof of \cite[Lemma 8]{zhao2024data} and our condition $\delta/\gamma \leq p_2$, it holds that $\|P-\hat{P}\| \leq p_1\delta/\gamma$. Since $\hat{B} \leq \|B\| + \delta/(\gamma)\leq \|B\| + p_2$, we have
		$
			\|B^{\top}P - \hat{B}^{\top}\hat{P}\|\leq (\|B\| + {p_2}){p_1\delta}/{\gamma} + {\|P\|\delta}/{\gamma}.
		$
		Further,  
		$\|E-\hat{E}\|
			\leq \|B\|\|P\|(1+\|K\|)\delta/\gamma+  ((\|B\| + {p_2})p_1 + {\|P\|} )
			   (\|A+BK\| +  (1+\|K\|){p_2} ){\delta}/{\gamma}.$
		
		Together with Lemma \ref{lem:bounds}, the proof is completed.
	\end{proof}
	
	Next, we quantify how the difference in the policy gradient in Lemma \ref{lem:grad_diff} affects the difference in the cost $|C(K'')-C(K')|$. We leverage the following preliminary lemma from \cite{fazel2018global}, which follows directly from Lemma \ref{lem:perturb}.
	\begin{lemma}\label{lem:cost_diff}
		Let $K \in \mathcal{S}$. There exist polynomials $p_4 = \text{poly}(C(K)/\underline{\sigma}(Q), \|A\|^{-1}, \|B\|^{-1}, \|R\|^{-1}, \underline{\sigma}(R) )$ and $p_5,p_6 = \text{poly}(C(K)/\underline{\sigma}(Q), \|A\|, \|B\|, \|R\|, 1/\underline{\sigma}(R))$ such that, if $\|\widetilde{K} - K\| \leq p_4$, then  $\widetilde{K} \in \mathcal{S}$ and
		$$
		\|\widetilde{\Sigma} - \Sigma\| \leq p_5\|\widetilde{K} - K\|, ~~
		|C(\widetilde{K}) - C(K)| \leq p_6\|\widetilde{K} - K\|.
		$$ 
	\end{lemma}
	
	Now, we bound the difference in the cost functions.
	 
	\begin{lemma}\label{lem:bias}
		Let $K \in \mathcal{S}$. There exists a polynomial $p_7$ in $({\underline{\sigma}(Q)}/{C(K)}, \|A\|^{-1},  \|B\|^{-1},  \|R\|^{-1}, \underline{\sigma}(R))$ such that, if
		$$
		\frac{\delta}{\gamma} \leq  p_2 ~~\text{and}~~ \eta \leq \min\left\{\frac{p_4\gamma}{p_3\delta}, p_7\right\},
		$$
		then it holds that
		$
		\left| C(K'') - C(K')\right| \leq  \eta p_3p_6 \delta/{\gamma}.
		$
	\end{lemma}
	\begin{proof}
		By our hypothesis and Lemma \ref{lem:grad_diff}, it holds that
		$$
		\|K' - K''\| = \eta \| \nabla{C}(K) - \nabla\hat{C}(K)\|\leq \frac{\eta p_3 \delta}{\gamma} \leq p_4.
		$$
		
		Let $p_7$ be the polynomial of the stepsize in the gradient descent case of \cite[Theorem 7]{fazel2018global}. Then, it follows from \cite[Theorem 7]{fazel2018global} that the update \eqref{equ:pg_truee} returns a stabilizing policy $K''$. Furthermore, we can apply Lemma \ref{lem:cost_diff} to obtain that
		$
		|C(K') - C(K'')| \leq p_6 \|K' - K''\| \leq {\eta p_3p_6 \delta}/{\gamma},
		$
		which completes the proof.
	\end{proof}
	
	With Lemma \ref{lem:bias} and \eqref{equ:exact}, we show the convergence of \eqref{equ:pg_appr}
	\begin{equation}\label{equ:appro}
		C(K') - C(K) \leq  -\frac{\eta}{2\mu} (C(K) - C^*) + \frac{\eta p_3p_6 \delta}{\gamma}.  
	\end{equation} 
	To  show the convergence of Algorithm \ref{alg:indirect}, we need to first prove that $C(K_t)$ is uniformly upper-bounded, such that the polynomials $p_i, i\in\{1,2,\dots, 7\}$ have uniformly bounds. Let  
	\begin{equation}\label{equ:bound_cost}
		\overline{C} = C^* + C(K_{t_0}) + 1 + \frac{1}{2\underline{l}\mu},
	\end{equation}
	where $\underline{l} := l(C^*)$, $\mu = {\|\Sigma^*\|}/{\underline{\sigma}(R)}$. Since the quantities $l$ and $p_i$ are function of $C(K)$,  let $\bar{l}, \bar{p}_1, \underline{p}_2, \bar{p}_3, \bar{p}_5, \bar{p}_6, \underline{p}_7$ be the associated quantities at $\overline{C}$, and $\underline{p}_4$ be the   quantity at $C^*$.  
	
	\begin{lemma}[Boundedness of the cost]\label{lem:bounded}
		If
		\begin{equation}\label{equ:cond}
			\frac{\delta_t}{\gamma_t} \leq  \min\left\{\underline{p}_2,  \frac{1}{2\bar{p}_6 \bar{p}_3\mu}\right\} ~ \text{and}~  \eta \leq \min\left\{\frac{\underline{p}_4\gamma_t}{\bar{p}_3\delta_t}, \underline{p}_7, \frac{1}{\bar{l}} \right\},
		\end{equation}
		then $C(K_t)$ of Algorithm \ref{alg:indirect} has a uniform upper bound, i.e.,
		$C(K_t) \leq \overline{C}$ with $\overline{C}$ defined in \eqref{equ:bound_cost}.
	\end{lemma} 
	\begin{proof}
		The proof is based on mathematical induction. Clearly, the bound holds at $t = t_0$, i.e., $C(K_{t_0}) \leq \overline{C}$. Suppose that $C(K_t) \leq \overline{C}$ for $t>t_0$. Next, we show $C(K_{t+1}) \leq \overline{C}$.
		
		By Lemma \ref{lem:bias}, \eqref{equ:exact}, and our hypothesis $C(K_t) \leq \overline{C}$, the gradient descent \eqref{equ:mbpg} yields
		\begin{align*}
			C(K_{t+1}) - C(K_t) &\leq -\frac{\eta}{2\mu} (C(K_t) - C^*) + \eta \bar{p}_6\bar{p}_3 \frac{\delta_{t+1}}{\gamma_{t+1}} \\
			&\leq -\frac{\eta}{2\mu} (C(K_t) - C^*) + \frac{\eta}{2\mu},
		\end{align*}
		where the last inequality follows from our condition on $\delta_t/\gamma_t$.
		
		Consider two cases. If $C(K_t) \geq C^* +1$, then 
		$$
		C(K_{t+1}) \leq C(K_t) - \frac{\eta}{2\mu}  + \frac{\eta}{2\mu} = C(K_t) \leq \overline{C}.
		$$
		Otherwise, if $C(K_t) < C^* +1$, then
		$$
		C(K_{t+1}) \leq C^* + 1 + \frac{\eta}{2\mu} \leq C^* + 1 + \frac{1}{2\bar{l}\mu}\leq \overline{C}.
		$$
		The proof is completed.
	\end{proof}
	

	Then, under the condition \eqref{equ:cond}, it follows that the convergence certificate \eqref{equ:conv} in Theorem \ref{thm:indirect} holds with $\nu_5 = \bar{p}_3\bar{p}_6$.  
	
	Next, we show the sequential stability of the closed-loop system of Algorithm \ref{alg:indirect}. To this end, we first find proper system-theoretic matrices for  $(\kappa, \alpha)$-strong stability.
	
	\begin{lemma}\label{lem:strong_stable}
		The policy $K \in \mathcal{S}$ is $(\kappa, \alpha)$-strongly stable with
		$$
		\kappa = \sqrt{\frac{C(K)}{\min\{\underline{\sigma}(R), \underline{\sigma}(Q)\}}},~\alpha = 1 - \sqrt{1-\frac{1}{\kappa^2}}.
		$$
	\end{lemma}
	
	\begin{proof}
		Let $H = \Sigma^{1/2}$ and $L = \Sigma^{-1/2} (A+BK) \Sigma^{1/2}$. Then, the closed-loop matrix satisfies $A+BK = HLH^{-1}$. 
		
		By the definition of $\Sigma$, it follows that 
		$
		I = \Sigma^{-1} + LL^{\top}.
		$
		Since $C(K) \geq \underline{\sigma}(Q) \text{Tr}(\Sigma)$, the following bounds hold
		$$
		\|L\| \leq \sqrt{1 - \frac{\underline{\sigma}(Q)}{C(K)}}, ~\|\Sigma^{\frac{1}{2}}\|\|\Sigma^{-\frac{1}{2}}\| \leq \sqrt{\frac{C(K)}{\underline{\sigma}(Q)}}.
		$$
		
		Noting $\|K\|\leq \sqrt{C(K)/\underline{\sigma}(R)}$, the proof is completed.
	\end{proof}
	
	We prove that the policy sequence of Algorithm \ref{alg:indirect}
	is sequentially strong stable.

	\begin{lemma}\label{lem:seq_stable}
		There exist $\bar{p}_{8}$ as a function of $\overline{C}$ such that, if
		\begin{equation}\label{equ:condition}
			\begin{aligned}
				&\frac{\delta_t}{\gamma_t} \leq  \min\left\{\underline{p}_2, \frac{1}{2\bar{p}_6 \bar{p}_3\mu}\right\} \\
				&\eta \leq \min\left\{\frac{\underline{p}_4\gamma_t}{\bar{p}_3\delta_t}, \underline{p}_7, \frac{1}{\bar{l}}, \frac{\underline{p}_4}{\bar{p}_8}, \frac{\underline{\alpha}}{4\bar{\kappa}
					^2 \bar{p}_5\bar{p}_8}, \frac{1}{2\bar{p}_5\bar{p}_8} \right\},
			\end{aligned}
		\end{equation} 
		then  $\{K_t\}$ of Algorithm \ref{alg:indirect} is sequentially strong stable with parameters $(\bar{\kappa}, \underline{\alpha})$, where $\bar{\kappa}, \underline{\alpha}$ are the quantities of Lemma \ref{lem:strong_stable} evaluated	 at $\overline{C}$.
	\end{lemma}
	
	\begin{proof}
		By Lemma \ref{lem:bounded}, the cost is uniformly bounded, i.e., $C(K_t)\leq \overline{C}$ for all $t\geq t_0$. Then, with the parameters $\bar{\kappa}, \underline{\alpha}$, the first two conditions (i) and (ii) in Definition \ref{def:sss} are satisfied, i.e., the policy $K_t$ is $(\bar{\kappa}, \underline{\alpha})$-strongly stable. Further, it suffices to show (iii) $\|H_{t+1}^{-1}H_t\| \leq 1 + \alpha/2$, or equivalently,
		$\|\Sigma_{t+1}^{-1}\Sigma_t\| \leq (1 + {\alpha}/{2})^2,$ where $\Sigma_t$ is the solution to \eqref{equ:Sigma} with $K_t$.
		
		By the perturbation theory for matrix inverse \cite[Theorem 35]{fazel2018global}, if $\|\Sigma_{t+1} - \Sigma_t\|<1/2$, then
		$
		\|\Sigma_{t+1}^{-1} - \Sigma_t^{-1}\| \leq {2\|\Sigma_{t+1} - \Sigma_t\|}/{\underline{\sigma}(\Sigma_t)} \leq 2\|\Sigma_{t+1} - \Sigma_t\|.
		$ Further,
		\begin{align*}
			\|\Sigma_{t+1}^{-1}\Sigma_t\|  
			& = \|(\Sigma_{t+1}^{-1} - \Sigma_t^{-1})\Sigma_t + I_n \|\\
			&\leq 1 + 2\|\Sigma_{t+1} - \Sigma_t\|\|\Sigma_t\| \\
			&\leq 1 + 2 \|\Sigma_{t+1} - \Sigma_t\|  {C(K_t)}/{\underline{\sigma}(Q)}  \\
			&\leq 1 + 2\bar{\kappa}^2 \|\Sigma_{t+1} - \Sigma_t\|.
		\end{align*}
		
		Thus, we require $ \|\Sigma_{t+1} - \Sigma_t\|$ to be sufficiently small. By Lemma 3, if $\|K_{t+1} - K_t\| \leq \underline{p}_4$, then $\|\Sigma_{t+1} - \Sigma_t\| \leq \bar{p}_5 \|K_{t+1} - K_t\|$. Thus, we need $\|K_{t+1} - K_t\| = \eta \|\nabla \hat{C}(K_{t+1})\|$ to be small. To this end, we provide a bound for $\|\nabla \hat{C}(K_{t+1})\|$. By Lemma \ref{lem:grad_diff}, we have 
		$
		\| \nabla\hat{C}(K) \| \leq \|\nabla{C}(K) \| + {\bar{p}_3\delta}/{\gamma} \leq \|\nabla{C}(K) \| + \bar{p}_3\underline{p}_2.
		$
		Since $C(K)\leq \overline{C}$, the right-hand side has a uniform upper bounded denoted by $\bar{p}_8$.
		Thus, as long as $\eta \leq \underline{p}_4/\bar{p}_8$, we have $\|K_{t+1} - K_t\| \leq \underline{p}_4$. To ensure $\|\Sigma_{t+1} - \Sigma_t\|<1/2$, we need $\eta \leq 1/(2\bar{p}_5\bar{p}_8)$. Furthermore, it follows that $\|\Sigma_{t+1}^{-1}\Sigma_t\| \leq  1 + 2\bar{\kappa}^2 \|\Sigma_{t+1} - \Sigma_t\| \leq 1 + 2\bar{\kappa}^2\bar{p}_5\bar{p}_8\eta$. Since we also require $2\kappa^2p_5p_8\eta \leq \alpha/2$, we let $\eta \leq \alpha/(4\kappa^2p_5p_8)$. Combining all the bounds on $\eta$ completes the proof.
	\end{proof}
	
	With sequential strong stability, we establish the boundedness of the state under Algorithm \ref{alg:indirect}. The proof follows the argument in \cite[Lemma 3]{cohen2019learning} and is included here for completeness.
	\begin{lemma}\label{lem:bound_state}
		Under the condition \eqref{equ:condition}, it holds that
		\begin{equation}\label{equ:state_bound}
			\|x_t\| \leq   \bar{\kappa}(1-\underline{\alpha}/2)^{t-t_0}\|x_{t_0}\|+\frac{2  \bar{\kappa}}{\underline{\alpha}} \max_{t_0\leq i<t}\|Be_i + w_i\|.
		\end{equation} 
	\end{lemma}
	 
	\begin{proof}
		Following the control policy $u_t = K_t x_t + e_t$, the state sequence satisfies
		$x_{t+1} = (A+BK_t)x_t + Be_t + w_t$, and 
		\begin{equation}\label{equ:state}
			x_t = F_{t_0}x_{t_0} + \sum_{i=t_0}^{t-1} F_{i+1}(Be_i + w_i),
		\end{equation} 
		where
		$
		F_i = \prod_{s = i}^{t-1}\left(A +B K_s\right), t_0\leq i \leq t-1,~\text{and}~ F_t = I_n.
		$
		
		By Lemma \ref{lem:seq_stable}, the policy sequence $\{K_t\}$ is $(\bar{\kappa}, \underline{\alpha})$-strongly stable. That is, there exist matrices $\{H_t\succ 0\}$ and $\{L_t\}$ such that $A+BK_t = H_tL_tH_t^{-1}$ for all $t\geq t_0$, and 
			(i) $\|L_t\|\leq 1 - \underline{\alpha},~\|K_t\|\leq \bar{\kappa}$;
			(ii) $\|H_t\|\leq \bar{\kappa}$, $\|H_t^{-1}\|\leq 1  $;
			(iii) $\|H_{t+1}^{-1}H_t\| \leq 1 + \underline{\alpha}/2. $
		Thus, for $t_0\leq i <t$, it holds that
		\begin{align*}
			\left\|F_i\right\| & =\left\|\prod_{s=i}^{t-1} H_s L_s^{\top} H_s^{-1}\right\| \\
			& \leq\left\|H_{t-1}\right\|\left(\prod_{s=i}^{t-1}\left\|L_s\right\|\right)\left(\prod_{s=i}^{t-2}\left\|H_{s+1}^{-1} H_s\right\|\right)\left\|H_i^{-1}\right\| \\
			& \leq \bar{\kappa}(1-\underline{\alpha})^{t-i}(1+\underline{\alpha}/2)^{t-i-1} \leq \bar{\kappa}(1-\underline{\alpha}/2)^{t-i}.
		\end{align*}
		
		Inserting it into \eqref{equ:state}, we obtain that
			\begin{align*}
				&\|x_t\|  \leq\|F_{t_0}\|\|x_{t_0}\|+\sum_{i={t_0}}^{t-1}\|F_{i+1}\|\|Be_i + w_i\| \\
				& \leq \bar{\kappa}(1-\frac{\underline{\alpha}}{2})^{t-t_0}\|x_{t_0}\|+ \bar{\kappa} \sum_{i=t_0}^{t-1}(1-\frac{\underline{\alpha}}{2})^{t-i-1}\|Be_i + w_i\| \\
				& \leq  \bar{\kappa}(1-\frac{\underline{\alpha}}{2})^{t-{t_0}}\|x_{t_0}\|+ \bar{\kappa} \max _{{t_0} \leq i<t}\|Be_i + w_i\| \sum_{t=1}^{\infty}(1-\frac{\underline{\alpha}}{2})^t \\
				& = \bar{\kappa}(1-\frac{\underline{\alpha}}{2})^{t-{t_0}}\|x_{t_0}\|+\frac{2  \bar{\kappa}}{\underline{\alpha}} \max_{{t_0} \leq i<t}\|Be_i + w_i\|.
			\end{align*}
		The proof is completed.
	\end{proof}
	
	Finally, we notice that the upper bound of $\eta$ in \eqref{equ:condition} depends on $\gamma_t/\delta_t$. By using the bound $\gamma_t/\delta_t \geq 1/\underline{p}_2$, it follows that $\eta$ has a constant upper bound, which completes the proof.
	
%
	
%
%
%
%
%
%
%

	\subsection{Proof of Lemma \ref{lem:equi}}
	By the covariance parameterization \eqref{equ:newpara} and the chain rule, the gradient satisfies $\Pi_{\overline{X}_{0}} \nabla J(V) = \Pi_{\overline{X}_{0}}\overline{U}_{0}^{\top}\nabla \hat{C}(K)$. Substituting it into \eqref{equ:pgdV} and using \eqref{equ:newpara} yield \eqref{equ:relation}. 
	
	To show that $\underline{\sigma}(M_t) \geq \gamma_t^4$, note that
	$$
	\Phi_t \Pi_{\overline{X}_{1,t}} \Phi_t^{\top} = \begin{bmatrix}
		M_t & 0_{m\times n} \\
		0_{n\times m} & 0_{n\times n}
	\end{bmatrix}.
	$$
	Then, it holds that $\underline{\sigma}(M_t) = \sigma_m (\Phi_t \Pi_{\overline{X}_{0,t}} \Phi_t^{\top})$. By the inequalities of singular value of matrix products, it follows that
	$
	\underline{\sigma}(M_t)
	\geq \underline{\sigma}(\Phi_t)\sigma_m(\Pi_{\overline{X}_{0,t}}\Phi_t^{\top})  
	\geq \underline{\sigma}(\Phi_t) \sigma_m(\Pi_{\overline{X}_{0,t}}) \underline{\sigma}(\Phi_t) 
	= (\underline{\sigma}(\Phi_t))^2\geq \gamma_t^2.
	$
	
	\subsection{Proof of Theorem \ref{thm:direct}} 
	We omit the subscript of $M_t$ when it is clear from the context. By Lemma \ref{lem:equi}, it suffices to consider the policy gradient update with estimated $(\hat{A}, \hat{B})$
	$$
		K' = K - \eta M \nabla\hat{C}(K)
	$$
	and the gradient update with the ground-truth $(A,B)$
	\begin{equation}\label{equ:pg_true}
		K'' = K - \eta M \nabla{C}(K).
	\end{equation}

We have the following convergence result for \eqref{equ:pg_true}.
\begin{lemma}\label{lem:pg}
	For $K\in \mathcal{S}$ and a stepsize $\eta \leq 1/(l\|M\|)$, the policy gradient update \eqref{equ:pg_true} satisfies
	$$
	C(K'') - C^* \leq \left( 1 - \frac{\eta \underline{\sigma}(M)}{2\mu}  \right)(C(K) - C^*).
	$$
\end{lemma}
\begin{proof}
	We start from a change of variables. Noting that $M$ is positive definite,  let $L = M^{-\frac{1}{2}}K$ and $\tilde{C}(L) := C(M^{\frac{1}{2}}L)$. Then, we have $\nabla \tilde{C}(L)=M^{\frac{1}{2}} \nabla C(K)$, and the update on $K$ in \eqref{equ:pg_true} is equivalent to the following update $
	L'' = L - \eta M^{\frac{1}{2}} \nabla C(K) = L - \eta \nabla \tilde{C}(L).$
	
	We show that the cost $\widetilde{C}(L)$ is gradient dominated and smooth in $L$. We have
		\begin{align*}
			& \left\|\nabla \tilde{C}(L)\right\|^2=\left\|M^{\frac{1}{2}} \nabla C(K)\right\|^2 \geq \underline{\sigma}(M)\left\|\nabla C(K)\right\|^2 \\
			&\geq   \frac{\underline\sigma(M)}{\mu}\left(C(K)-C^*\right) = \frac{\underline{\sigma}(M)}{\mu}\left(\tilde{C}(L)-C^*\right),
		\end{align*}
	where the second inequality follows from Lemma \ref{lem:pl}. Thus, the gradient dominance property holds with constant $\underline{\sigma}(M)/\mu$.
	
	Moreover, for any stabilizing $L_1,L_2$ satisfying $L+\psi(L'-L), \forall \psi \in [0,1]$ is stabilizing, it holds that
		\begin{align*}
			& \left\|\nabla \tilde{C}(L_1)-\nabla \tilde{C}\left(L_2\right)\right\|=\left\|M^{\frac{1}{2}} \nabla C(K_1)-M^{\frac{1}{2}} \nabla C\left(K_2\right)\right\| \\
			& \leq  l\|M\|^{\frac{1}{2}}\left\|K_1-K_2\right\| \leq l\|M\|\left\|L_1-L_2\right\|,
		\end{align*}
	that is, $\widetilde{C}(L)$ is smooth with smoothness constant $l\|M\|$. Noting $\widetilde{C}(L) = C(K)$ completes the proof.
\end{proof}

\begin{lemma}\label{lem:bias_deepo}
	Let $K \in \mathcal{S}$. If
	$$
	\frac{\delta}{\gamma} \leq   p_2  ~~\text{and}~~ \eta \leq \frac{1}{\|M\|}\cdot \min\left\{\frac{p_4\gamma}{p_3\delta}, {p_7}\right\},
	$$
	then it holds that
	$
	\left| C(K'') - C(K')\right| \leq {\eta p_6p_3 \delta \|M\|}/{\gamma}.	$
\end{lemma}
\begin{proof}
	By our hypothesis and Lemma \ref{lem:grad_diff}, it holds that
	$$
	\|K'' - K'\| \leq \eta \|M\| \| \nabla\hat{C}(K) - \nabla C(K)\|\leq \frac{\eta p_3 \delta\|M\|}{\gamma} \leq p_4.
	$$
	Following the same vein of the proof of Lemma \ref{lem:bias}, we can show that $K''$ is stabilizing, and
	$
	|C(K') - C(K'')| \leq p_6 \|K' - K''\| \leq {\eta p_6p_3 \delta \|M\|}/{\gamma},
	$
	which completes the proof.
\end{proof}

Next, we show that the cost of Algorithm \ref{alg:deepo} has an upper bound, i.e., $C(K_t) \leq \overline{C}$, where $\overline{C}$ is given by \eqref{equ:bound_cost}.

\begin{lemma}[Boundedness of the cost]\label{lem:bounded_direct}
	If
	\begin{equation}\label{equ:cond_direct}
		\begin{aligned}
		&\frac{\delta_t}{\gamma_t} \leq  \min\left\{\underline{p}_2,   \frac{\underline{\sigma}(M_{t})}{2\bar{p}_6 \bar{p}_3\mu\|M_{t}\|}\right\}  \\ &\eta_t \leq \frac{1}{\|M_{t}\|}\cdot \min\left\{\frac{\underline{p}_4\gamma_t}{\bar{p}_3\delta_t}, \underline{p}_7, \frac{1}{\bar{l}} \right\},
		\end{aligned} 
	\end{equation}
	then $C(K_t)$ has a uniform upper bound, i.e.,
	$C(K_t) \leq \overline{C}$.
\end{lemma} 

\begin{proof}
	The proof is based on mathematical induction. Clearly, the bound holds at $t = t_0$, i.e., $C(K_{t_0}) \leq \overline{C}$. Suppose that $C(K_{t}) \leq \overline{C}$ for $t > t_0$. Next, we show $C(K_{t+1}) \leq \overline{C}$.
	
	By Lemmas \ref{lem:equi}, \ref{lem:pg}, \ref{lem:bias_deepo} and our hypothesis on $\delta/\gamma$, we have
	\begin{align*}
		&C(K_{t+1}) - C(K_t) \\
		&\leq -\frac{\eta_{t+1}\underline{\sigma}(M_{t+1})}{2\mu} (C(K_t) - C^*) + \frac{\eta_{t+1}\bar{p}_6\bar{p}_3 \delta_{t+1}\|M_{t+1}\|}{\gamma_{t+1}} \\
		&\leq -\frac{\eta_{t+1}\underline{\sigma}(M_{t+1})}{2\mu} (C(K_t) - C^*) + \frac{\eta_{t+1}\underline{\sigma}(M_{t+1})}{2\mu}.
	\end{align*} 
	 
	Consider two cases. If $C(K_t) \geq C^* + 1$, then 
	$$
	C(K_{t+1}) \leq C(K_t) - \frac{\eta_{t+1}\underline{\sigma}(M_{t+1})}{2\mu}   + \frac{\eta_{t+1}\underline{\sigma}(M_{t+1})}{2\mu}  \leq \overline{C}.
	$$
	Otherwise, if $C(K_t) < C^* + 1$, then
	$$
	C(K_{t+1}) \leq C(K_t)  + \frac{\eta_{t+1}\underline{\sigma}(M_{t+1})}{2\mu}  \leq C^* + 1+ \frac{\underline{\sigma}(M_{t+1})}{2\bar{l}\mu\|M_{t+1}\|}.
	$$ 
	Noting that $\underline{\sigma}(M_{t+1}) \leq \|M_{t+1}\|$, it further leads to $C(K_{t+1})\leq \overline{C}$. By induction, the proof is completed.
\end{proof}
 
Furthermore, we  show the convergence of Algorithm \ref{alg:deepo}. Under the condition \eqref{equ:cond_direct},  we have
\begin{align*}
	&C(K_t) - C^* \\
	&\leq (1-\frac{\eta_t\underline{\sigma}(M_{t})}{2\mu}) (C(K_{t-1}) - C^*) + \frac{\eta_t \bar{p}_6\bar{p}_3\delta_t \|M_{t}\|}{\gamma_{t}}\\
	&\leq (1-\frac{\eta_t\underline{\sigma}(M_{t})}{2\mu}) (C(K_{t-1}) - C^*) + \frac{ \bar{p}_6\bar{p}_3\delta_t }{\bar{l} \gamma_t},
\end{align*}
where the last inequality follows from $\eta_t \leq 1/(\bar{l}\|M_{t}\|)$. Then, the convergence certificate in Theorem \ref{thm:direct} holds with $\nu_6 = {\bar{p}_6\bar{p}_3}/{\bar{l}}$.
Finally, we prove that the policy sequence $\{K_t\}$ of Algorithm \ref{alg:deepo}
is sequential strong stable.

\begin{lemma}\label{lem:seq_stable_direct}
	Suppose that 
		\begin{align*}
			&\frac{\delta_t}{\gamma_t} \leq  \min\left\{\underline{p}_2,   \frac{\underline{\sigma}(M_{t})}{2\bar{p}_6 \bar{p}_3\mu\|M_{t}\|}\right\} ~\text{and}\\
			&\eta_t \leq \frac{1}{\|M_{t}\|} \cdot \min\left\{\frac{\underline{p}_4\gamma_t}{\bar{p}_3\delta_t}, \underline{p}_7, \frac{1}{\bar{l}}, \frac{\underline{p}_4}{\bar{p}_8}, \frac{\underline{\alpha}}{4\bar{\kappa}
				^2 \bar{p}_5\bar{p}_8}, \frac{1}{2\bar{p}_5\bar{p}_8} \right\}.
		\end{align*} 
	then  $\{K_t\}$ of Algorithm \ref{alg:deepo} is sequentially strong stable with parameters $(\bar{\kappa}, \underline{\alpha})$. Moreover, the state is bounded as \eqref{equ:state_bound}.
\end{lemma}

\begin{proof}
The proof follows the same vein of that of Lemmas \ref{lem:seq_stable} and \ref{lem:bound_state}, and is omitted due to space limitation.
\end{proof}

 Noting that  the condition on stepsize in \eqref{lem:seq_stable} can be further simplified by using $\gamma_t/\delta_t \geq 1/\underline{p}_2$, the proof is completed.

\section{Proofs in Section \ref{sec:grad}}
\subsection{Proof of Lemma \ref{lem:bridge}}
 	 	First, we derive the natural policy gradient descent over $V$.  Let $\psi = \text{vec}(V^{\top})$. Then, the covariance parameterization between $K$ and $V$ \eqref{equ:newpara} can be expressed in terms of $\theta$ and $\psi$ 
\begin{equation}\label{equ:para}
	\theta   = \text{vec}\left(V^{\top} \overline{U}_0^{\top}\right) = \left(\overline{U}_0 \otimes I_n\right)\text{vec}(V^{\top}) = \left(\overline{U}_0\otimes I_n\right)\psi
\end{equation} 
with a subspace constraint 
	$(\overline{X}_0 \otimes I_n)\psi = \text{vec}(I_n). $
Thus, the FIM with respect to $\psi$ follows from the transformation rule
\begin{equation}\label{equ:trans}
	F_\psi =  \frac{\partial \theta^{\top}}{\partial \psi} F_\theta \frac{\partial \theta}{\partial \psi} 
	= \left(\overline{U}_0\otimes I_n\right)^{\top} F_{\theta} \left(\overline{U}_0 \otimes I_n\right).
\end{equation}
Since $\psi$ lies in an affine space, we should restrict the Fisher inverse to its tangent space $\mathcal{N}(\overline{X}_0\otimes I_n)$. Denote $O$ as a unit orthogonal basis of $\mathcal{N}(\overline{X}_0)$. Then, the corresponding orthogonal basis of $\mathcal{N}(\overline{X}_0\otimes I_n)$ is $(O\otimes I_n)$, and the Fisher inverse along the tangent space is given by 
$
((O\otimes I_n)^{\top}F_\psi (O\otimes I_n))^{-1}.$ Finally, the natural gradient method updates $\psi$ in the feasible affine space as $\psi^+ = \psi - \eta \Delta \psi$, where
$
\Delta \psi= (O\otimes I_n) ((O\otimes I_n)^{\top}F_\psi (O\otimes I_n))^{-1} (O\otimes I_n)^{\top} \nabla J(\psi).
$

Next, we show that the change in $\psi$ leads to the change of $\theta$ in \eqref{equ:ng}. By the parameterization \eqref{equ:para}, $\theta^+ =(\overline{U}_0\otimes I_n)\psi^+ = \theta - \eta (\overline{U}_0\otimes I_n) \Delta \psi$. By the chain rule, it holds that $\nabla J(\psi) = (\overline{U}_0\otimes I_n)^{\top} \nabla C(\theta)$. Note that by Lemma \ref{lem:equi}, the matrix $(\overline{U}_0\otimes I_n)(O\otimes I_n) = (\overline{U}_0O)\otimes I_n $ is invertible. Then, combining \eqref{equ:trans} and the expression of $\Delta \psi$, it follows that the update of $\theta$ is $\theta^+ = \theta -  \eta F_{\theta}^{-1}\nabla C(\theta)$, which is exactly \eqref{equ:ng}.

Finally, noting that the above derivations hold under the estimated model $(\hat{A}_{t+1}, \hat{B}_{t+1})$, the proof is completed.

\subsection{Proof of Theorem \ref{thm:na}}
The proof follows the same vein as that of Theorem \ref{thm:indirect}. 
Consider the policy gradient update with estimated $(\hat{A}, \hat{B})$
\begin{equation}\label{equ:pg_appro_na}
		K' =  K - 2\eta \hat{E}
\end{equation} 
and the update with the ground-truth $(A,B)$
\begin{equation}\label{equ:pg_true_na}
	K'' =   K - 2\eta {E},
\end{equation} 
where $E$ is defined in \eqref{equ:E}. We first quantify the difference between the exact and approximated natural gradient $\hat{E}-E$, and then quantify $C(K'')-C(K')$.

\begin{lemma}\label{lem:grad_diff_na}
	Let $K \in \mathcal{S}$. There exists a polynomial $p_9=\text{poly}(C(K)/\underline{\sigma}(Q), \|A\|, \|B\|, \|R\|, 1/\underline{\sigma}(R))$ such that if $\delta / \gamma \leq  p_2$, then
	$
		\| \hat{E}-E \| \leq {p_9\delta}/{\gamma}.
	$ 
\end{lemma}

The proof follows from that of Lemma \ref{lem:grad_diff} and is omitted.

\begin{lemma}\label{lem:bias_na}
	Let $K \in \mathcal{S}$ and $p_{10} = {1}/{\|R+B^{\top}PB\|}$. If
	$$
	\frac{\delta}{\gamma} \leq  p_2 ~~\text{and}~~ \eta \leq  \min\left\{\frac{p_4\gamma}{2p_9\delta},  p_{10}\right\},
	$$
	then it holds that
	$
	\left| C(K'') - C(K')\right| \leq {\eta p_6p_9 \delta}/{\gamma}.
	$
\end{lemma}
\begin{proof}
	By our hypothesis and Lemma \ref{lem:grad_diff_na}, it holds that
	$
	\|K' - K''\| = 2\eta \| \hat{E}-E\|\leq {2\eta p_9 \delta}/{\gamma} \leq p_4.
	$
	By taking $\eta \leq p_{10}$, it follows from \cite[Theorem 7]{fazel2018global} that the update \eqref{equ:pg_true_na} returns a stabilizing policy $K''$. Furthermore, we can apply Lemma \ref{lem:cost_diff} to obtain that
	$
	|C(K') - C(K'')| \leq p_6 \|K' - K''\| \leq {\eta p_6p_9 \delta}/{\gamma},
	$
	which completes the proof.
\end{proof}

By \cite[Lemma 15]{fazel2018global}, if $\eta \leq p_{10}$, then the progress of the exact natural gradient descent \eqref{equ:pg_true_na} satisfies
$$
	C(K'') - C(K) \leq  -\frac{\eta}{2\mu} (C(K) - C^*).  
$$
Together with Lemma \ref{lem:bias_na}, the progress of \eqref{equ:pg_appro_na} is
\begin{equation}\label{equ:cost_na_appro}
	C(K'') - C(K) \leq  -\frac{\eta}{2\mu} (C(K) - C^*) + \frac{\eta p_6p_9 \delta}{\gamma}.  
\end{equation} 
Define $
	\overline{C} = C^* + C(K_0) + 1 + \frac{1}{2\|R\|\mu}$. Then, we show that the cost is uniformly upper bounded by $\overline{C}$.

\begin{lemma}[Boundedness of the cost]\label{lem:bounded_na}
If
\begin{equation}\label{equ:cond_na}
	\frac{\delta_t}{\gamma_t} \leq  \min\left\{\underline{p}_2,  \frac{1}{2\bar{p}_6 \bar{p}_9\mu}\right\} ~\text{and}~ \eta \leq \min\left\{\frac{\underline{p}_4\gamma_t}{2\bar{p}_9\delta_t}, \underline{p}_{10} \right\},
\end{equation}
then $C(K_t)$ of PGAC with natural gradient has a uniform upper bound, i.e.,
$C(K_t) \leq \overline{C}, \forall t\geq t_0$.
\end{lemma} 
\begin{proof}
	The proof is based on mathematical induction. Clearly, the bound holds at $t = t_0$, i.e., $C(K_{t_0}) \leq \overline{C}$. Suppose that $C(K_{t}) \leq \overline{C}$ for $t > t_0$. Next, we show $C(K_{t+1}) \leq \overline{C}$.
	
	By \eqref{equ:cost_na_appro} and our hypothesis $C(K_t) \leq \overline{C}$, the gradient descent \eqref{equ:ngK_indirect} yields
	$
		C(K_{t+1}) - C(K_t)  \leq -\frac{\eta}{2\mu} (C(K_t) - C^*) + \eta \bar{p}_6\bar{p}_9 \frac{\delta_{t+1}}{\gamma_{t+1}} 
		\leq -\frac{\eta}{2\mu} (C(K_t) - C^*) + \frac{\eta}{2\mu},
	$
	where the last inequality follows from our condition on $\delta_t/\gamma_t$.
	
	Consider two cases. If $C(K_t) \geq C^* +1$, then 
	$$
	C(K_{t+1}) \leq C(K_t) - \frac{\eta}{2\mu}  + \frac{\eta}{2\mu} = C(K_t) \leq \overline{C}.
	$$
	Otherwise, if $C(K_t) < C^* +1$, then
	$$
	C(K_{t+1}) \leq C^* + 1 + \frac{\eta}{2\mu} \leq C^* + 1 + \frac{1}{2\|R\|\mu}\leq \overline{C},
	$$
	where the second inequality follows from $\eta \leq \underline{p}_{10} \leq 1/\|R\|$.
	The proof is completed.
\end{proof}
 
Then, under the condition \eqref{equ:cond_na}, it holds that
	\begin{align*}
		C(K_{t}) - C^* &\leq (1-\frac{\eta}{2\mu})^{t-t_0} (C(K_{t_0}) - C^*) \\
		&+ \eta \bar{p}_6\bar{p}_9\sum_{i = t_0}^{t-1}(1-\frac{\eta}{2\mu})^{t-1-i} \frac{1}{\text{SNR}_i}.
	\end{align*} 

\begin{lemma}\label{lem:seq_stable_na}
	There exists $\bar{p}_{11}$ as a function of $\overline{C}$ such that, if
	\begin{equation}\label{equ:condition_na}
		\begin{aligned}
			&\frac{\delta_t}{\gamma_t} \leq \min\left\{\underline{p}_2,   \frac{1}{2\bar{p}_6 \bar{p}_9\mu}\right\},~ \text{and}\\
			&\eta \leq \min\left\{\frac{\underline{p}_4\gamma_t}{2\bar{p}_9\delta_t}, \underline{p}_{10}, \frac{\underline{p}_4}{2\bar{p}_{11}}, \frac{\underline{\alpha}}{8\bar{\kappa}
				^2 \bar{p}_5\bar{p}_{11}}, \frac{1}{4\bar{p}_5\bar{p}_{11}} \right\},
		\end{aligned}
	\end{equation}
	then  $\{K_t\}$   is sequentially strong stable with parameters $(\bar{\kappa}, \underline{\alpha})$, where $\bar{\kappa}, \underline{\alpha}$ are the quantities at $\overline{C}$. Moreover, the state is bounded as \eqref{equ:state_bound}.
\end{lemma}

\begin{proof}
  	The proof follows the same vein as that of Lemma \ref{lem:seq_stable}, where we require $ \|\Sigma_{t+1} - \Sigma_t\|$ to be sufficiently small. By Lemma 3, if $\|K_{t+1} - K_t\| \leq \underline{p}_4$, then $\|\Sigma_{t+1} - \Sigma_t\| \leq \bar{p}_5 \|K_{t+1} - K_t\|$. Thus, we need $\|K_{t+1} - K_t\| = 2\eta \|\hat{E}_{t+1}\|$ to be small. To this end, we provide a bound for $\|E_{t+1}\|$. By Lemma \ref{lem:grad_diff_na}, we have 
	$
	\|\hat{E}_{t+1} \| \leq \|E_{t+1}\| + {\bar{p}_9\delta_t}/{\gamma_t} \leq \|E_{t+1}\| + \bar{p}_9\underline{p}_2.
	$
	Since $C(K_t)\leq \overline{C}$, the right-hand side has a uniform upper bounded denoted by $\bar{p}_{11}$.
	Thus, as long as $\eta \leq \underline{p}_4/(2\bar{p}_{11})$, we have $\|K_{t+1} - K_t\| \leq \underline{p}_4$. To ensure $\|\Sigma_{t+1} - \Sigma_t\|<1/2$, we need $\eta \leq 1/(4\bar{p}_5\bar{p}_{11})$. Furthermore, it follows that $\|\Sigma_{t+1}^{-1}\Sigma_t\| \leq  1 + 2\bar{\kappa}^2 \|\Sigma_{t+1} - \Sigma_t\| \leq 1 + 2\bar{\kappa}^2 \bar{p}_5 \|K_{t+1}-K_t\|  \leq 1 + 4\bar{\kappa}^2\bar{p}_5\bar{p}_{11}\eta$. Since we also require $4\kappa^2p_5p_{11}\eta \leq \alpha/2$, we let $\eta \leq \alpha/(8\kappa^2p_5p_{11})$. Combining all the bounds on $\eta$ completes the proof.
\end{proof}

\subsection{Proof of Theorem \ref{thm:indirect_gn} }
Let $\eta = 1/2$. Consider the Gauss-Newton update with estimated $(\hat{A}, \hat{B})$
\begin{equation}\label{equ:pg_appr_gn}
	K' = K - (R+\hat{B}^{\top}\hat{P}\hat{B})^{-1}\hat{E},
\end{equation} 
and the update with the ground-truth $(A,B)$ 
\begin{equation}\label{equ:pg_true_gn}
	K'' = K- (R+B^{\top}PB)^{-1}E. 
\end{equation}
where $E$ is defined in \eqref{equ:E}. We first quantify the difference between the exact and approximated Gauss-Newton gradient, and then quantify $C(K'')-C(K')$.

\begin{lemma}\label{lem:grad_diff_gn}
	 Let $K \in \mathcal{S}$. Then, there exist $p_{12}$ as a monotonously decreasing function of $C(K)$ and
	$p_{13}=\text{poly}(C(K)/\underline{\sigma}(Q), \|A\|, \|B\|, \|R\|, 1/\underline{\sigma}(R))$, such that if $\delta / \gamma \leq  p_{12}$, then $
		\| (R+\hat{B}^{\top}\hat{P}\hat{B})^{-1}\hat{E}-(R+B^{\top}PB)^{-1}E \| \leq {p_{13}\delta}/{\gamma}.$
\end{lemma}

The proof follows from that of Lemma \ref{lem:grad_diff} and is omitted.

\begin{lemma}\label{lem:bias_gn}
	Let $K \in \mathcal{S}$. If
	$
	{\delta}/{\gamma} \leq  \min\{p_{12}, {p_4}/{p_{13}} \},
	$
	then it holds that
	$
	| C(K'') - C(K')| \leq { p_6p_{13} \delta}/{\gamma}.
	$
\end{lemma}
\begin{proof}
	By our hypothesis and Lemma \ref{lem:grad_diff_gn}, it holds that
	$
		\|K' - K''\|  = \| (R+\hat{B}^{\top}\hat{P}\hat{B})^{-1}\hat{E}-\left(R+B^{\top}PB\right)^{-1}E \| \leq {p_{13} \delta}/{\gamma} \leq p_4.
	$
	Then, we can apply Lemma \ref{lem:cost_diff} to obtain that
	$
	|C(K') - C(K'')| \leq p_6 \|K' - K''\| \leq {  p_6p_{13}\delta}/{\gamma},
	$
	which completes the proof.
\end{proof}

By \cite[Lemma 8]{fazel2018global}, the Gauss-Newton method \eqref{equ:pg_true_gn} yields
$$
	C(K'') - C(K) \leq  -\frac{\eta}{\|\Sigma^*\|} (C(K) - C^*).  
$$
Together with Lemma \ref{lem:bias_gn}, the progress of \eqref{equ:pg_appr_gn} is
\begin{equation}\label{equ:cost_gn}
	C(K') - C(K) \leq  -\frac{\eta}{\|\Sigma^*\|} (C(K) - C^*) + \frac{  p_6p_{13}\delta}{\gamma}.  
\end{equation}

Let  $
	\overline{C} = C^* + c + \frac{c}{2\|\Sigma^*\|}$ for any constant $c>0$. Then, we have the following result.
\begin{lemma}[Boundedness of the cost]\label{lem:bounded_gn}
	If
	\begin{equation}\label{equ:cond_gn}
	\frac{\delta_t}{\gamma_t} \leq  \min\left\{\underline{p}_{12},  \frac{\underline{p}_4}{\bar{p}_{13}}, \frac{c}{2\bar{p}_6\bar{p}_{13}\|\Sigma^*\|} \right\}, ~ C(K_{t_0} ) \leq \overline{C},
	\end{equation}
	then $C(K_t)$ of Algorithm \ref{alg:Hewer} satisfies
	$C(K_t) \leq \overline{C}$.
\end{lemma} 
\begin{proof}
	The proof is based on mathematical induction. By our condition on $C(K_0)$, the bound holds at $t = 0$. Suppose that $C(K_t) \leq \overline{C}$. Next, we show $C(K_{t+1}) \leq \overline{C}$.
	
	By \eqref{equ:cost_gn}  and our hypothesis $C(K_t) \leq \overline{C}$, the Gauss-Newton method \eqref{equ:GNK_indirect} yields
	\begin{align*}
		C(K_{t+1}) - C(K_t) &\leq -\frac{1}{2\|\Sigma^*\|} (C(K_t) - C^*) +  \frac{  \bar{p}_6\bar{p}_{13}\delta_{t+1}}{\gamma_{t+1}} \\
		&\leq -\frac{1}{2\|\Sigma^*\|} (C(K_t) - C^*) + \frac{c}{2\|\Sigma^*\|} ,
	\end{align*}
	where the last inequality follows from our condition on $\delta_t/\gamma_t$.
	
	Consider two cases. If $C(K_t) \geq C^* +c$, then 
	$
	C(K_{t+1}) \leq C(K_t) - {\eta}/{(2\mu)}  + {\eta}/{(2\mu)} = C(K_t) \leq \overline{C}.
	$
	Otherwise, if $C(K_t) < C^* +c$, then
	$
	C(K_{t+1}) \leq C^* + c + {c}/{(2\|\Sigma^*\|)} = \overline{C}.
	$
	The proof is completed.
\end{proof}

\begin{lemma}\label{lem:seq_stable_gn}
 There exists a constant $c$ and $\bar{p}_{13}$ as a function of $\overline{C}$ such that, if
	\begin{equation}\label{equ:condition_gn}
	\frac{\delta_t}{\gamma_t} \leq  \min\left\{\underline{p}_{12},   \frac{\underline{p}_4}{\bar{p}_{13}}, \frac{c}{2\bar{p}_6\bar{p}_{13}\|\Sigma^*\|}, \frac{c}{\bar{p}_{13}} \right\},
	\end{equation} 
	then  $\{K_t\}$ of Algorithm \ref{alg:Hewer} is sequentially strong stable with parameters $(\bar{\kappa}, \underline{\alpha})$, where $\bar{\kappa}, \underline{\alpha}$ are the quantities at $\overline{C}$.
\end{lemma}

\begin{proof}
	The proof follows the same vein as that of Lemma \ref{lem:seq_stable}, where we require $ \|\Sigma_{t+1} - \Sigma_t\|$ to be sufficiently small. By Lemma 3, if $\|K_{t+1} - K_t\| \leq \underline{p}_4$, then $\|\Sigma_{t+1} - \Sigma_t\| \leq \bar{p}_5 \|K_{t+1} - K_t\|$. Thus, we need $\|K_{t+1} - K_t\| =  \|(R+\hat{B}_{t+1}^{\top}\hat{P}_{t+1}\hat{B}_{t+1})^{-1}\hat{E}_{t+1}\|$ to be small. To this end, we provide a bound for $\|(R+\hat{B}^{\top}\hat{P}\hat{B})^{-1}\hat{E}\|$. By Lemma \ref{lem:grad_diff_gn}, we have 
	$
	\|(R+\hat{B}^{\top}\hat{P}\hat{B})^{-1}\hat{E} \| \leq \|(R+B^{\top}PB)^{-1}E\| + {\bar{p}_{13}\delta}/{\gamma}.
	$
	
	By \cite[Lemma 11]{fazel2018global}, it holds that
	$
	\|(R+B^{\top}PB)^{-1}E\| \leq  (   {(C(K)- C^*)}/{\underline{\sigma}(R)}  )^{\frac{1}{2}}.
	$
	Since $C(K_t)\leq \overline{C}$, we have
	\begin{equation}\label{equ:ccc}
			\begin{aligned}
			\|K_{t+1} - K_t\| &= \|(R+\hat{B}_{t+1}^{\top}\hat{P}_{t+1}\hat{B}_{t+1})^{-1}\hat{E}_{t+1}\| \\
			&\leq \left( \frac{c}{\underline{\sigma}(R)} + \frac{c}{2\|\Sigma^*\|\underline{\sigma}(R)} \right)^{\frac{1}{2}} + c.
		\end{aligned}
	\end{equation}

	Thus, as long as the right-hand side of the this inequality is smaller than $\underline{p}_4$, we have $\|K_{t+1} - K_t\| \leq \underline{p}_4$. To ensure $\|\Sigma_{t+1} - \Sigma_t\|<1/2$, we additionally require $\|K_{t+1} - K_t\| \leq 1/(2\bar{p}_5)$. Since $\|\Sigma_{t+1}^{-1}\Sigma_t\| \leq  1 + 2\bar{\kappa}^2 \|\Sigma_{t+1} - \Sigma_t\| \leq 1 + 2\bar{\kappa}^2 \bar{p}_5 \|K_{t+1}-K_t\|$, to ensure sequential stability we need $2\bar{\kappa}^2 \bar{p}_5 \|K_{t+1}-K_t\|\leq \underline{\alpha}/2$. Noting that $\|K_{t+1} - K_t\|$ is upper bounded as \eqref{equ:ccc}, these conditions can be ensured for a sufficiently small $c$.
\end{proof}

The rest of the proof follows the same vein as that of Theorem \ref{thm:indirect} and is omitted.

\section{Proof of Theorem \ref{thm:indirect_reg}}
The idea is to show that, under the given selection of $\lambda_t$, the difference between the policy gradient $\nabla{C}(K)$ and the regularized one $\nabla\hat{C}_{\lambda}(K)$ scales linearly with $\delta/\gamma$. Hence, the results follow the proof of Theorem \ref{thm:indirect}.
\begin{lemma}\label{lem:grad_diff_reg}
	Let $K \in \mathcal{S}$. Then, there exist $p_{14}=\text{poly}(C(K)/\underline{\sigma}(Q), \|A\|, \|B\|, \|R\|, 1/\underline{\sigma}(R))$ and a constant $c_0>0$ such that, if $\delta / \gamma \leq  p_2$ and $\lambda \leq  c_0\delta$, then
	$
	\| \nabla\hat{C}_{\lambda}(K) - \nabla{C}(K) \| \leq {p_{14}\delta}/{\gamma}.
	$
\end{lemma}
\begin{proof}
	The difference in the gradient can be bounded by
	\begin{equation}\label{equ:deco}
		\begin{aligned}
			&\| \nabla\hat{C}_{\lambda}(K) - \nabla{C}(K) \| \\
			&\leq \| \nabla\hat{C}_{\lambda}(K) -  \nabla {C}_{\lambda}(K) \| + \|  \nabla {C}_{\lambda}(K) - \nabla{C}(K) \|.
		\end{aligned}
	\end{equation} 
	
	We first bound the second term of \eqref{equ:deco}. By Lemma \ref{lem:grad_reg}, it holds
	$
	\left\| \nabla {C}_{\lambda}(K) - \nabla {C}(K)\right\| 
	= 2\lambda\|(\Phi^{-1})_{uu}K + B^{\top}\tilde{P}(A+BK) +  (\Phi^{-1})_{ux}\|\|\Sigma\|,
	$
	where $\tilde{P}$ is the solution to
	$
	\tilde{P} = (\Phi^{-1})_{xx} + K^{\top}(\Phi^{-1})_{uu}K + K^{\top}(\Phi^{-1})_{ux} + (\Phi^{-1})_{xu}K 
	+ (A+BK)^{\top}\tilde{P}(A+BK).
	$ For the first term, we notice that $\|(\Phi^{-1})_{uu}\|\leq \|\Phi^{-1}\|$ and $\|(\Phi^{-1})_{ux}\|\leq \|\Phi^{-1}\|$, and
	$
	\|\tilde{P}\| \leq \|\Sigma\|\|(\Phi^{-1})_{xx} + K^{\top}(\Phi^{-1})_{uu}K + K^{\top}(\Phi^{-1})_{ux} + (\Phi^{-1})_{xu}K\| \leq  \|\Phi^{-1}\| \|\Sigma\| (1+ \|K\|^2 + 2\|K\|).
	$
	Then, there exists $p_{15}=\text{poly}(C(K)/\underline{\sigma}(Q), \|A\|, \|B\|, \|R\|, 1/\underline{\sigma}(R))$ such that
	$
	\left\| \nabla {C}_{\lambda}(K) - \nabla {C}(K)\right\|
	\leq p_{15} \lambda \|\Phi^{-1}\| \leq p_{15}c_0 {\delta}/{\gamma},
	$
	where  the last inequality follows from $\|\Phi^{-1}\| = 1/\underline{\sigma}(\Phi) \leq 1/\gamma$ and our condition on $\lambda$.
	 Analogously, we can show that an upper bound of $ \| \nabla\hat{C}_{\lambda}(K) -  \nabla {C}_{\lambda}(K)  \|$ is linear in $\delta/\gamma$. Inserting these bounds in \eqref{equ:deco} completes the proof. 
\end{proof}

The rest of the proof follows the same vein as that of Theorem \ref{thm:indirect} and is omitted due to space limitation.

	%

\bibliographystyle{IEEEtran}
\bibliography{mybibfile}

\begin{IEEEbiography}[{\includegraphics[width=1in,height=1.25in,clip,keepaspectratio]{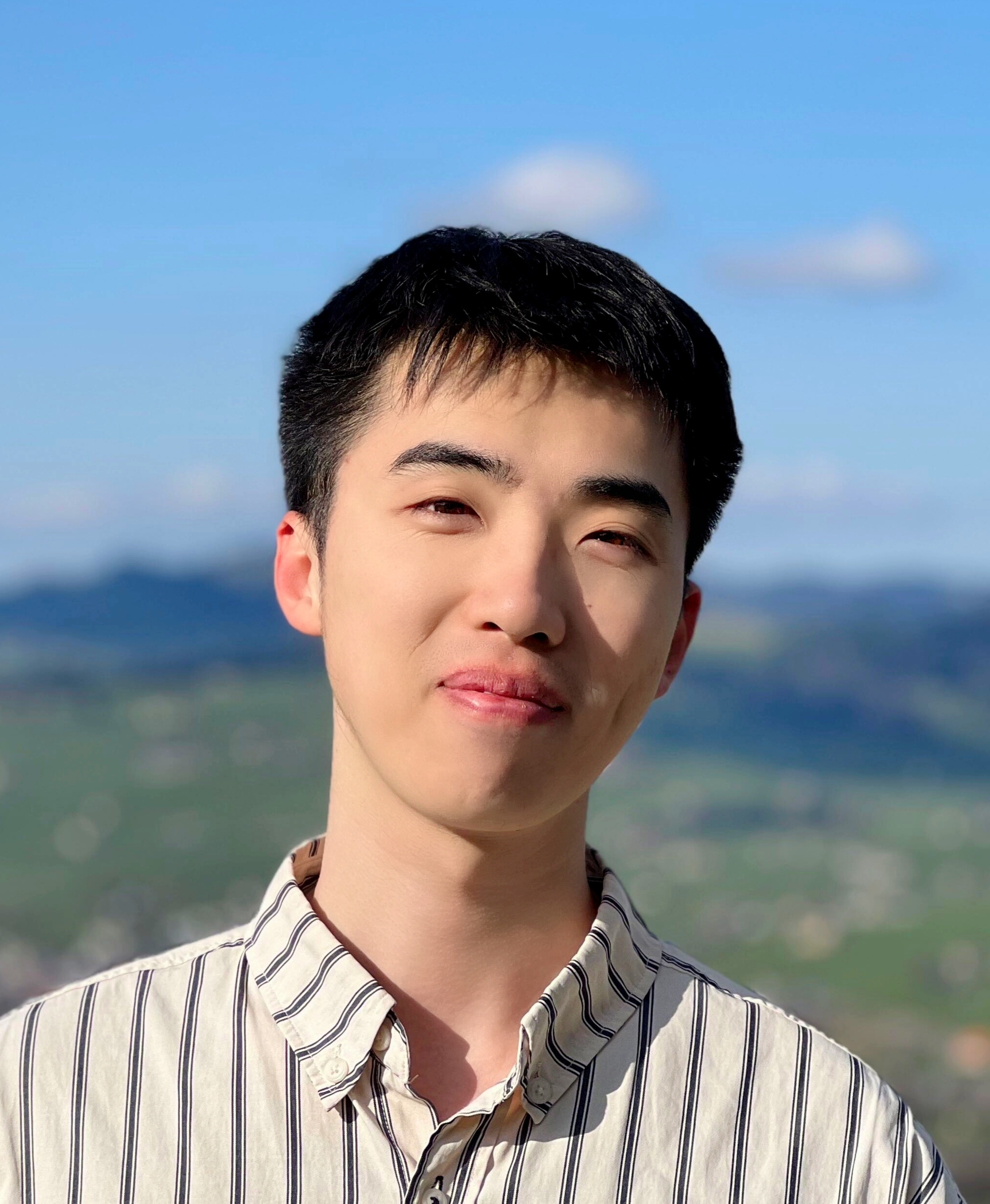}}]{Feiran Zhao} received the B.S. degree in  Control Science and Engineering from the Harbin Institute of Technology, China, in 2018, and the Ph.D. degree in Control Science and Engineering from the Tsinghua University, China, in 2024. He is currently a postdoctoral researcher at ETH Z\"{u}rich. His research interests include policy optimization, data-driven control, adaptive control and their applications.
\end{IEEEbiography}

\begin{IEEEbiography}
	[{\includegraphics[width=1in,height=1.25in,clip,keepaspectratio]{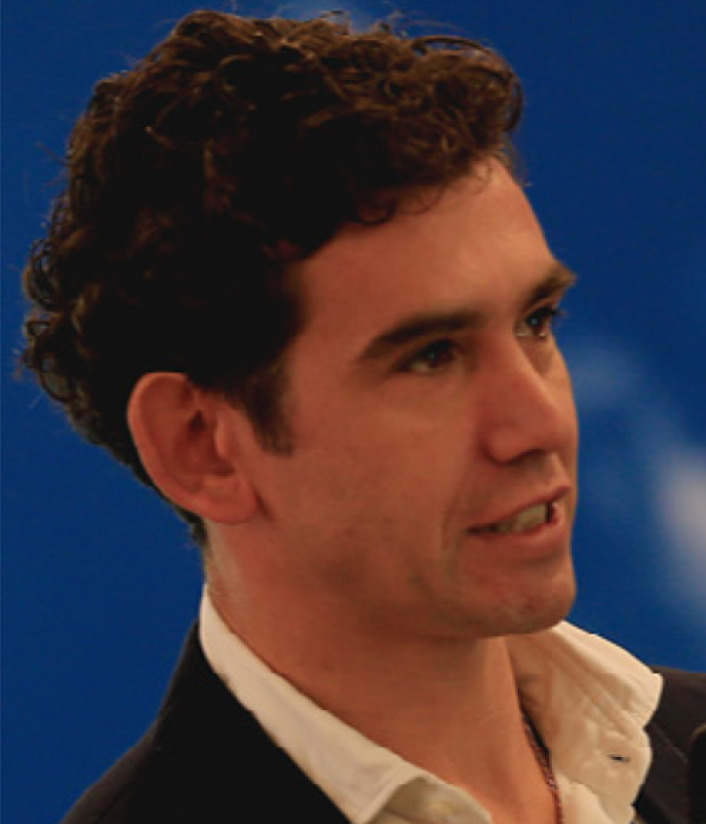}}]
	{Alessandro Chiuso} (Fellow, IEEE) is Professor with the Department of Information Engineering, Universit\`a di Padova. He received the ``Laurea" degree summa cum laude in Telecommunication Engineering from the University of Padova in July 1996 and the Ph.D. degree (Dottorato di ricerca) in System Engineering from the University of Bologna in 2000. He has been a visiting scholar with the Dept. of Electrical Engineering, Washington University St. Louis and Post-Doctoral fellow with the Dept. Mathematics, Royal Institute of Technology, Sweden. He joined the University of Padova as an Assistant Professor in 2001, Associate Professor in 2006 and then Full Professor since 2017. He currently serves as Editor (System Identification and Filtering) for Automatica. He has served as an Associate Editor for Automatica, IEEE Transactions on Automatic Control, IEEE Transactions on Control Systems Technology, European Journal of Control and MCSS. He is chair of the IFAC Coordinating Committee on Signals and Systems. He has been General Chair of the IFAC Symposium on System Identification, 2021 and he is a Fellow of IEEE (Class 2022). His research interests are mainly at the intersection of Machine Learning, Estimation, Identification and their Applications, Computer Vision and Networked Estimation and Control.
\end{IEEEbiography}

\begin{IEEEbiography}
	[{\includegraphics[width=1in,height=1.25in,clip,keepaspectratio]{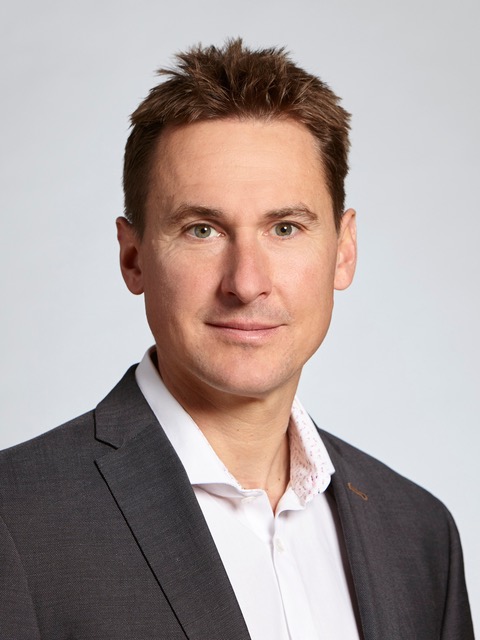}}]
			{Florian D\"{o}rfler} is a Professor at the Automatic Control Laboratory at ETH Z\"{u}rich. He received his Ph.D. degree in Mechanical Engineering from the University of California at Santa Barbara in 2013, and a Diplom degree in Engineering Cybernetics from the University of Stuttgart in 2008. From 2013 to 2014 he was an Assistant Professor at the University of California Los Angeles. He has been serving as the Associate Head of the ETH Z\"{u}rich Department of Information Technology and Electrical Engineering from 2021 until 2022. His research interests are centered around automatic control, system theory, and optimization. His particular foci are on network systems, data-driven settings, and applications to power systems. He is a recipient of the distinguished young research awards by IFAC (Manfred Thoma Medal 2020) and EUCA (European Control Award 2020). His team has received many best thesis and best paper awards in the top venues of control, power systems, power electronics, and circuits and systems. He is currently serving on the council of the European Control Association and as a senior editor of Automatica.
	\end{IEEEbiography}	

\end{document}